\input amstex.tex
\documentstyle{amsppt}
 
\input graphicx.tex

\magnification=1200
\hsize=150truemm
\vsize=224.4truemm
\hoffset=4.8truemm
\voffset=12truemm

\TagsOnRight
\NoBlackBoxes
\NoRunningHeads

\def\Square{\rlap{$\sqcup$}$\sqcap$}
\def\cqfd {\quad \hglue 7pt\par\vskip-\baselineskip\vskip-\parskip
{\rightline{\Square}}}

\define\Z{{\bold Z}}
\define\N{{\bold N}}
\define\Rt {${\bold R}$-tree}
\define\mi{^{-1}}
\define\st   {\text{\rm Stab}\,}
\define\rk{\text{\rm rk}\,}
\define\fix{\text{\rm Fix}\,}
\let\thm\proclaim
\let\fthm\endproclaim
\let\inc\subset 
\let\ds\displaystyle
\let\ev\emptyset

\let\ov\overline

\define\bo{\partial }

\define\Aut{\text{\rm Aut}\,}
\define\Out{\text{\rm Out}\,}

\define\cc#1{\ov{#1}}
\define\fd{f_\sharp}

\newcount\tagno
\newcount\secno
\newcount\subsecno
\newcount\stno
\global\subsecno=1
\global\tagno=0
\define\ntag{\global
\advance\tagno by 1\tag{\the\tagno}}

\newcount\figno
\newcount\fihno
\global\figno=0
\global\fihno=1
\define\fig{\global\advance\figno by 1 \global\advance\fihno by 1
{\the\figno}}

\define\sta{\uppercase\expandafter{
\the\secno}.\the\stno
\global\advance\stno by 1}

\define\sect{\global\advance\secno by 1\global\subsecno=1\global\stno=1\
\uppercase\expandafter{
\the\secno}. }

\def\nom#1{\edef#1{\uppercase\expandafter{
\the\secno}.\the\stno}}
\def\eqnom#1{\edef#1{(\the\tagno)}}

\newcount\refno
\global\refno=0
\def\nextref#1{\global\advance\refno by 1\xdef#1{\the\refno}}
\def\bref {\ref\global\advance\refno by 1\key{\the\refno}}


 \nextref\BFH
\nextref\BFHd
\nextref\BG
\nextref\BH
\nextref\CT
\nextref\DS
\nextref\GJLL
\nextref\GL
\nextref\GLL
\nextref\GaLu
\nextref\Ln
\nextref\LLD
\nextref\Pig
\nextref\Sh

\topmatter

\title  Counting growth types of automorphisms of free groups
   \endtitle

\author  Gilbert Levitt 
\endauthor

\abstract Given an automorphism of a free group $F_n$, we consider the
following invariants: $e$ is the number of exponential strata (an upper bound for
the number of different exponential growth rates of conjugacy classes); $d$ is
the maximal degree of polynomial growth of conjugacy classes; $R$ is the rank of
the fixed subgroup. We determine precisely which triples $(e,d,R)$ may be
realized by an automorphism of $F_n$.  In particular, the inequality
$e\le\frac{3n-2}4$ (due to Levitt-Lustig) always holds. In an appendix, we show
that any conjugacy class grows like a polynomial times an exponential under
iteration of the automorphism.
\endabstract

\endtopmatter

\document \secno=0\stno=1

\head Introduction\endhead

Consider an automorphism $\alpha $ of a free group $F_n$ which is induced by a
homeomorphism $h$ of a compact surface $\Sigma $. After replacing $h $ by
a power, it becomes isotopic to  a homeomorphism $h'$ with a very simple
structure:
$\Sigma $ is a union of invariant subsurfaces $\Sigma _i$, and on each subsurface
$h'$ is either pseudo-Anosov, or a Dehn twist in an annulus, or the identity.  

Three very different behaviors thus appear:  exponential, linear, trivial. This
qualitative analysis also has a quantitative side.  Just writing that the sum of the
Euler characteristics of the $\Sigma _i$'s is $\chi  (\Sigma )$ gives bounds for,
say, the rank of the fixed subgroup of $\alpha $, or the number of subsurfaces. 

If $\alpha $ is an arbitrary automorphism of $F_n$, there is no general analogue
of the Nielsen-Thurston decomposition into invariant subsurfaces. Still, one can
again distinguish three different behaviors. Exponential behavior comes from
exponentially growing strata in a relative train track representative of $\alpha
$, or, more intrinsically, from the attracting laminations. Linear should now be
replaced by polynomial: unlike in the surface case, it is   possible for the
length of a conjugacy class to grow as a polynomial of degree $>1$ under
iteration of $\alpha
$ (the simplest example is the automorphism of $F_3$ which sends $a$ to $a$, $b$
to $ba$, $c$ to
$cb$, with $c$ growing quadratically). Trivial behavior comes from the fixed
subgroup of
$\alpha
$. 

Though  these three different behaviors are not separated as clearly as in the surface case,
the goal of this paper is to show that one can still give
precise numerical bounds for invariants that measure how much ``space'' each of
these three behaviors occupies within $F_n$. 

The first invariant attached to an automorphism $\alpha $ is the number $e$ of
attracting laminations (equal to the number of exponential strata of an
improved train track representative). It is an upper bound for the number of
exponential growth rates of conjugacy classes (see appendix). In the surface
case,
$e$ is the number of subsurfaces $\Sigma _i$ on which the map is pseudo-Anosov
(it is a pleasant exercise to compute the maximal value of $e$ on a given
orientable surface
$\Sigma
$; the answer is given in Remark 5.1).

The second invariant (which has no equivalent in the surface case) is the maximal
degree $d$ such that the length of some conjugacy class  grows as a
polynomial of degree $d$ under iteration of $\alpha $ (it equals 2 in the example
given above on $F_3$).

The third invariant is the rank $R=\rk\fix\alpha $ of the fixed subgroup of
$\alpha $. It is a famous theorem by Bestvina-Handel [\BH] that $R$ does not
exceed 
$n$.

Our main result is a precise characterization of which triples $(e,d,R)$ may be
realized by an automorphism of $F_n$. Let us first consider $e$ and $d$.

\thm{Theorem 1} Given $\alpha \in Aut(F_n)$, the numbers $e$ and $d$ satisfy:
$$\aligned
e+d&\le n-1\\
4e+2d&\le 3n-2\qquad\text{($\le 3n-3$ if $d>0$)}. 
\endaligned$$
In particular, $\ds e\le\frac{3n-2}4$ (Levitt-Lustig, [\Ln]).

Conversely, any $(e,d)$ satisfying these inequalities may be realized by some
$\alpha \in\Aut(F_n) $.
\fthm

The inequalities are equivalent to saying  that $(e,d)$ belongs to the closed
quadrilateral with vertices   $(0,0),(0,n-1),(\frac{n-1}2,\frac{n-1}2),
(\frac{3n-2}4,0)$ pictured on Figure \fig. 
The equality $e=
\left[\frac{3n-2}4\right]$  is achieved by a surface
homeomorphism with all subsurfaces $\Sigma _i$
either once-punctured tori or four-punctured spheres (see Figure \fig{}  for a
picture when $n\equiv2\mod 4$); recall that there is no pseudo-Anosov map on a
pair of pants.

\thm{Theorem 2} Given   $e$ and $d$ satisfying the conditions above, the
possible values of $\rk\fix\alpha $ for an automorphism $\alpha $ of $F_n$ are
exactly those allowed by the following inequalities:
$$\aligned
e+\max(d-1,0)+\rk\fix\alpha &\le n  \\
4e+2d+2\rk\fix\alpha &\le 3n+1\qquad \text{($\le 3n$ if $d=0$)}.
\endaligned$$
\fthm

\midinsert 
\centerline 
{\includegraphics[scale=.6]
{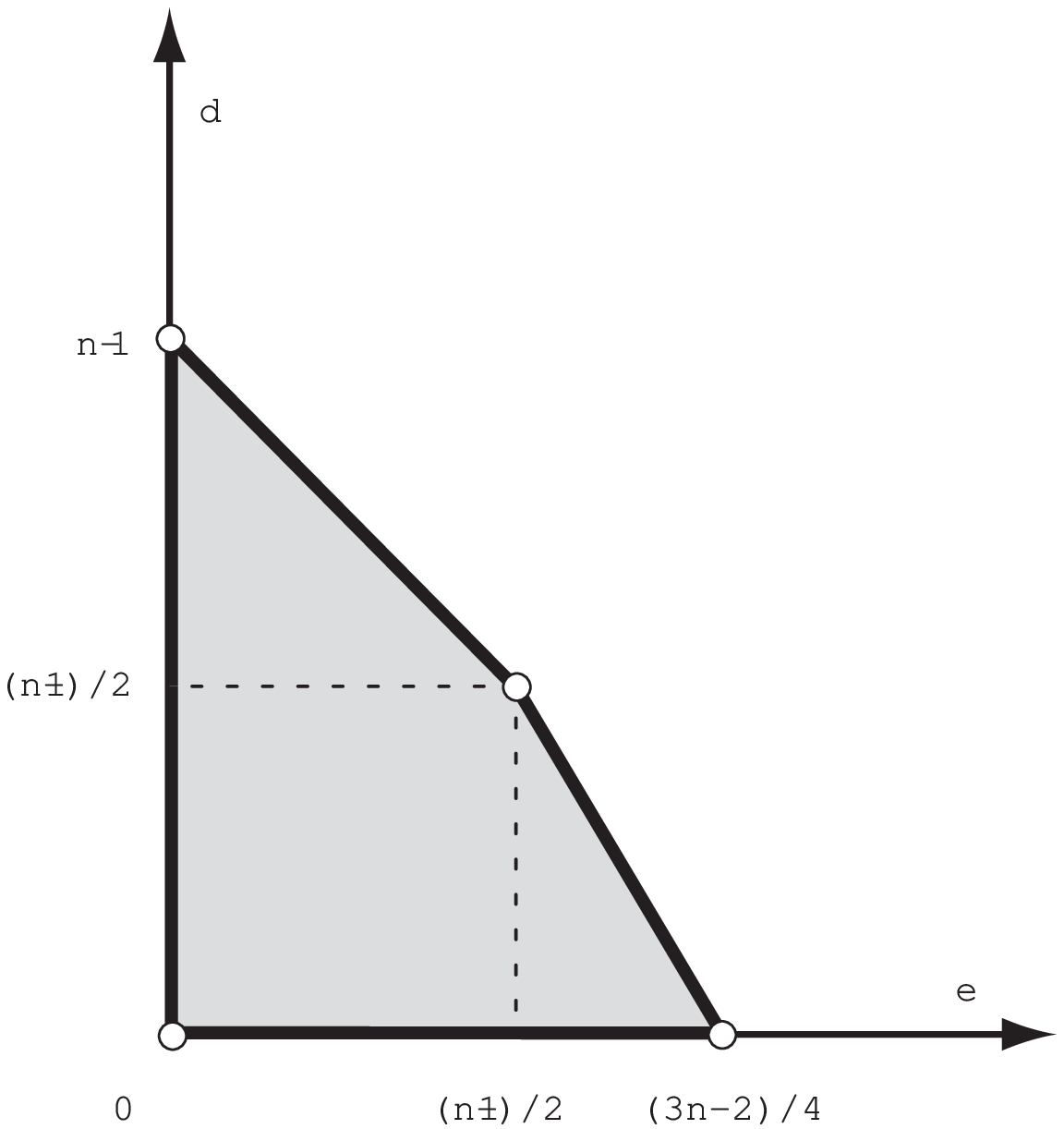}} 
\captionwidth{220pt}
\botcaption 
 {Figure 1}{ Possible values for $(e,d)$.}
\endcaption
\endinsert 

\midinsert 
\centerline 
{\includegraphics[scale=.5]  
{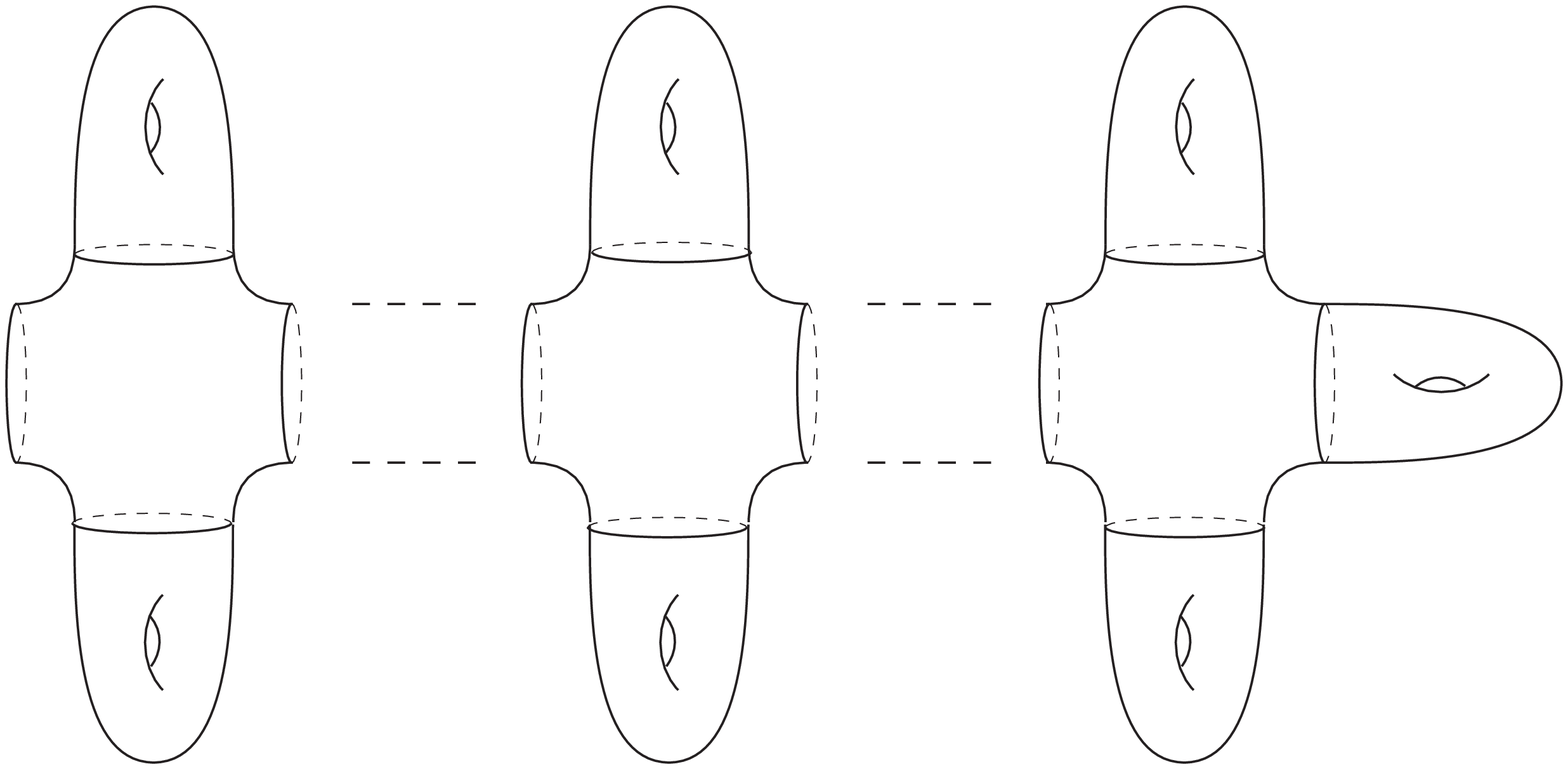}}
\captionwidth{230pt}
\botcaption 
 {Figure 2}{A surface $\Sigma $ with $\pi _1(\Sigma )\simeq F_n$ 
  decomposed into $\frac{3n-2}4$ subsurfaces carrying pseudo-Anosov
maps.}
\endcaption
\endinsert

The first inequality is a strengthening of the bound $\rk\fix\alpha \le n$ proved
in [\BH]. If the rank equals $n$, one must have $e=0$ and $d\le1$, as proved
in [\CT].

Our results   have to do with growth. In an appendix, we show:

\thm {Theorem 3}  Given
$\alpha \in\Aut(F_n)$ and an  element (or conjugacy class) $g$ in $F_n$, 
there exist $\lambda \ge1$ and $m\in\N$, and $C_1,C_2>0$, such that 
the length
of
$\alpha ^p(g)$  is bounded between $C_1\lambda ^pp^m$ and $C_2\lambda ^pp^m$
for all $p\ge1$.
\fthm

This result is
not deep (if one assumes train tracks), but has never appeared in print. Our
proof uses arguments from    a preliminary version of [\BG]. We also explain how to determine
the growth types $(\lambda ,m)$ with $\lambda >1$   from
the set of   attracting laminations of $\alpha $ and the Perron-Frobenius
eigenvalues. As shown in [\LLD], the numbers $\lambda $ may be viewed as
H\"older exponents associated to periodic points on the boundary of $F_n$.

Given $\alpha $, the number $e'$ of growth
types
$(\lambda ,m)$ with $\lambda >1$   is bounded by $e$, and {\it  
Theorems 1 and 2 are valid with $e$ replaced by $e'$.}
We also show (Theorem 4.7) that  the degree $m$ is bounded by $n/2-1$ when
$\lambda >1$ (when
$
\lambda =1$, the optimal bound for $m$ is $  n-1$ for growth of conjugacy
classes, $n$ for growth of elements of $F_n$).

Let us now say a few words about the proofs. As mentioned above, one of our
results is a strengthening of the ``Scott conjecture'' $\rk\fix\alpha \le n$.
As in [\GLL],  we work with the outer automorphism $\Phi   $
determined by $\alpha $, and instead of $\rk\fix\alpha $ we use a related
invariant
$r$ associated to $\Phi $.  In the surface case, $r$ equals the absolute value of the
Euler characteristic of the whole subsurface $\Sigma _f$ where $h'$ equals the
identity, whereas
$\rk\fix\alpha $ can only see one   component of $\Sigma _f$.

The proof is by induction on $n$. If
$\Phi  $ is polynomially growing (i.e\. $e=0$), some power of $\Phi  $ preserves
(up to conjugacy) a nontrivial decomposition of $F_n$ as a free product, or
preserves a free factor of rank $n-1$. This makes induction possible. 

If $e>0$, we consider a 
$\Phi $-invariant \Rt{} $T$ with trivial arc stabilizers, as in [\GLL]. Point
stabilizers have rank $<n$   and we can argue by induction.  The induction
starts from  the polynomial subgroups of $\Phi $:
a canonical finite family of conjugacy classes of subgroups on which all the
polynomial growth is concentrated    [\GaLu]. Our techniques  give
precise bounds for the ranks of these subgroups (see Theorem 4.1). 

Going back to the tree $T$, the inequality that would make the induction work to
prove our main results is $\sum_{i=1}^b (3n_i-2)\le 3n-6$, where $G_1,\dots, 
G_b$ are   representatives for   conjugacy classes of non-trivial  stabilizers
of points of $T$, and $n_i=\rk G_i$. Unfortunately, this inequality is false, as the
following example shows.

Suppose that $\Phi $ is induced by a homotopy equivalence $f$ of a finite complex
$Y$ obtained by attaching the boundary of  a once-punctured  torus $S$ onto a
graph
$\Gamma
$ with
$\pi _1\Gamma \simeq F_{n-1}$. Assume furthermore that $S$ and $\Gamma
$ are $f$-invariant, and $f_{|S}$  is   pseudo-Anosov.
There is a $\Phi $-invariant tree $T$ in which $\pi
_1\Gamma $ is a stabilizer, so that $\sum  (3n_i-2)=3(n-1)-2=3n-5$  
(there is a similar example  with $S$ a four-punctured sphere). 

To make
the induction work,   we need to have a more precise control on  the
outer automorphism $\Psi $ induced on  $\pi _1\Gamma$
  than
  stated so far. 
The key remark is the following. The boundary of $S$ provides a
   nontrivial conjugacy class $\cc\gamma $ which is fixed by $\Psi $.
 The element $\gamma $ is a commutator in
$F_n$, but it cannot be one in $\pi _1\Gamma $. Indeed, it must be a primitive
element of $\pi _1\Gamma $, as otherwise   $\pi _1Y$ could not be free. This
motivates the introduction of an invariant $k$, computed from periodic
conjugacy classes, so that   an inequality  
$\sum_i(3n_i-2-k_i)\le 3n-6-k$   does hold (Proposition 3.1).

In the first, preliminary, section, we explain how the induction works, and we
define the polynomial subgroups. In Section 2, we define the basic invariants and
we study how they behave under the induction process. In Section 3, we establish
the needed inequalities about the numbers $n_i$.
We prove the   inequalities of our main theorems in Section 4, as Corollaries 4.5
and 4.6. In Section 5  we construct examples, showing that the inequalities are
optimal. This section may be read independently of the others. As mentioned
above, we study growth in an appendix.

\vskip8pt
{\baselineskip=3pt
{\eightpoint Acknowledgement. {\it This research started in 1996 as joint work
with Martin Lustig, and the inequality
$e\le\frac{3n-2}4$ is proved in [\Ln]. Though Martin refused to coauthor the
present paper, he certainly deserves a lot of credit for it.

I also thank the referee for a suggestion which simplified the proof of Proposition 3.1 a great deal.
}}\par}

\head  \sect Preliminaries \endhead

 We view $\Phi \in\Out(F_n)$ as a collection of automorphisms  $\alpha
\in\Aut(F_n)$. We say that $\alpha $ {\it represents\/}
$\Phi
$, and we write $\Phi =\hat\alpha $. We write $i_x$ for the inner automorphism
$g\mapsto xgx\mi$, so that $\hat\beta =\hat\alpha $ if and only if $\beta
=i_x\circ \alpha $ for some $x\in F_n$. We say that $\alpha $ and $\beta $ are {\it
isogredient\/} if $\beta =i_y\circ\alpha \circ (i_y)\mi$ for some $y$ (the word
``similar'' was used in [\GJLL] and [\GLL]). 

We write $\fix\alpha $ for the fixed subgroup $\fix\alpha =\{g\in F_n\mid \alpha
(g)=g\}$, and $\rk\fix\alpha $ for its rank. Isogredient automorphisms have
conjugate  fixed subgroups.

\subhead Train tracks and laminations [\BFH]\endsubhead

Given $\Phi \in\Out(F_n)$, there exists $q\ge1$ such that $\Phi ^q$ is represented
by an {\it improved
relative train track map\/}    $f:G\to G$, as in [\BFH, Theorem 5.1.5].
We   denote by
$H_i$ the {\it $i$-th stratum}. The {\it height\/} of a path $\gamma $ is the largest
$i$ such  that
$\gamma $ contains an edge of $H_i$. The image of an edge of $H_i$ is a path of
height $\le i$. 

There are three types of strata. 
If $e $ is an edge in a 0-stratum $H_i$, then $f(e )$ has height $<i$. If $H_i$ is
an NEG stratum, it consists of a single edge $e_i$, and $f(e_i)=e_iu_i$ with $u_i$
of height $<i$. If $H_i$ is exponential,  it has a {\it transition matrix\/} $M$ whose
entry
$M_{pq}$ records the number of times that the image of the $p$-th edge of $H_i$
crosses the $q$-th edge (in either direction). This matrix is positive, and it has a
largest eigenvalue $\lambda >1$ called the {\it Perron-Frobenius eigenvalue\/} of
$H_i$. It is an algebraic integer.

To an exponential stratum $H_i$ is associated an {\it attracting lamination\/}
$\Lambda _i$ of $\Phi $ [\BFH, Section 3]. It may be defined as the set of
bi-infinite paths
$\gamma
$ in
$G$ such that any finite subpath of $\gamma $ is contained in some
tightened image $\fd^p(e)$, for $p\ge1$ and $e$ an edge of $H_i$. These paths
$\gamma $ are the {\it leaves\/} of $\Lambda _i$. Through the identification
between
$\pi _1(G)$ and $F_n$, one may view $\Lambda _i$ as a lamination on $F_n$ (a
subset of the quotient of $\bo F_n\times \bo F_n$ minus the
diagonal by  
 the action of $\Z/2\Z$, which interchanges the factors, and the diagonal
action of $F_n$). 

We have described the attracting laminations through the choice of a train
track representative $f$, but the set of attracting laminations on $F_n$
depends only on
$\Phi
$. For any $f$ representing a power of $\Phi $, there is a bijection between the set
of attracting laminations of
$\Phi
$ and the set of exponential strata of $f$.

\subhead Growth \endsubhead

We denote by $\cc g$ the {\it conjugacy class\/} of $g\in F_n$. The outer
automorphism
$\Phi $ acts on the set of conjugacy classes of elements, and on the set
of conjugacy classes of subgroups. 

If we fix a free basis for $F_n$,
the {\it length} $|g|$ is the length of the reduced word representing $g$. The
length
$|\cc g|$   is the length of any cyclically reduced word representing a conjugate of
$g$.  

We say that two sequences $a_p$ and $b_p$ are equivalent, or that $a_p$
{\it grows like\/} $b_p$, if the ratio $a_p/b_p$ is bounded away from $0$ and
infinity. It is well-known that any sequence $|\alpha ^p(g)|$ (resp\. $|\Phi  ^p(\cc
g)|$) grows like a polynomial of degree $m\in \N$, or has exponential growth. We
then say that
$g$ (resp\.  $\cc g$) grows {\it polynomially with degree $m$, or exponentially\/},
under iteration of $\alpha $ (resp\.  $\Phi $).  This does not depend on the choice
of a   basis. 

If
$g$ belongs to a finitely generated 
$\alpha
$-invariant   subgroup $H$, the growth of $g$ under $\alpha _{|H}$ is the same as
its growth under $\alpha $, because $H$ is quasiconvex in $F_n$ [\Sh]. Similarly,
the growth of
$\cc g$ under
$\Phi  _{|H}$ is the same as its growth under $\Phi  $.

An automorphism $\Phi $ (or a representative $\alpha $) is {\it polynomially
growing\/} if every conjugacy class grows polynomially. If $f$ represents a
power of $\Phi $, this is equivalent to saying that $f$ has no exponential stratum.

These facts are sufficient for the proof of Theorems 1 and 2. In the
appendix, we shall prove the more precise result that $g$ and $\cc g$   always
grow like a sequence
$\lambda ^pp^m$ with $\lambda \ge1$ and $m\in\N$. We say that $(\lambda ,m)$ is the {\it growth
type\/}. If $\lambda =1$, the growth
is polynomial. If
$\lambda >1$, the growth is   exponential and we say that $\lambda $ is the
{\it exponential growth rate\/}.  If
$\Phi ^q
$ is represented by an improved train track map, then any exponential growth
rate is the $q$-th root of the Perron-Frobenius eigenvalue of an exponential
stratum.

\subhead Setting up the induction \endsubhead

We shall prove our main results by induction on $n$. First we replace $\Phi
\in\Out(F_n)$ by a power,  so that it is represented by an improved relative
train track map $f:G\to G$.

First suppose that  $\Phi $ is polynomially growing. We consider the highest
stratum. It   consists of a single edge $E$. 

 \nom\refer
\thm{Definition \sta}
Let $\Phi $ be polynomially growing, represented by $f:G\to G$.
\roster
\item
If there is 
 a non-trivial
decomposition $F_n=G_1*G_2$ such that $\Phi $ has a representative  
$\alpha $ which
 satisfies $\alpha (G_i)=G_i$, we let $n_i$ be the rank of $G_i$ and we
denote $\alpha _i=\alpha _{|G_i}$, so that $\alpha =\alpha _1*\alpha _2$.
This happens in particular whenever the top edge $E$ separates $G$.
\item
If $E$ does not separate, 
there is a decomposition $F_n=G_1*\langle t\rangle$ (with $G_1$ of rank
$n-1$) and a representative $\alpha $ of $\Phi $ such that $\alpha (G_1)=G_1$
and $\alpha (t)=tu$ with $u\in G_1$.  We
  assume that $u$ cannot be written $u=a\alpha (a \mi)$ with $a\in G_1$, since
otherwise $\alpha (ta)=ta$ and we reduce to the previous  case $\alpha =\alpha
_1*\alpha _2$. We denote
$\alpha _1=\alpha _{|G_1}$. 
 \endroster 
\fthm

This will allow us 
to deduce properties of $\alpha $ from
properties of the $\alpha _i$'s.

In the remainder of this subsection, we assume that
  $\Phi $ is not polynomially growing. We then use 
 the $\Phi $-invariant \Rt{} $T$ associated to the highest exponential
stratum $H$, in the following sense.

If
$H$ is the highest of all strata, $T$ is   the tree constructed in
section 2 of [\GJLL]. If   not, we adapt the construction of [\GJLL] as follows. Let
$H'$ be the union of 
$H$ and all (non-exponential) strata above it. We consider the 
transition matrix
$M_{H'}$: there is one row and one column for each edge of $H'$,  and each
entry  records the number of times that the image of an edge  crosses an
edge (in either direction).
The eigenvalues
of  $M_{H'}$   are those of the transition matrix $M_H$ of $H$,
together with 0's and 1's. In particular, $M_{H'}$ has a  non-negative
eigenvector associated to the Perron-Frobenius eigenvalue $\lambda >1$
of $M_H$. We then apply the construction of [\GJLL], using this eigenvector
to give a PF-length to edges of $H'$ (edges below $H$ have PF-length 0).

 In either case, the group
$F_n$ acts on the \Rt{} $T$ isometrically, minimally, with no global fixed point.
All
  arc stabilizers are trivial.

 We let $\ell:F_n\to[0,\infty)$ be the length
function of
$T$, defined by $\ell(g)=\min_{x\in T}d(x,gx)$. Its value on $g$
only depends on the conjugacy class 
$\cc g$, so we sometimes write   $\ell(\cc g)$.
An element
$g$, or its class $\cc g$, is elliptic if $\ell(g)=0$, hyperbolic if
$\ell(g)>0$. An elliptic element has a unique fixed point, a hyperbolic
element has an axis along which it acts as translation by $\ell(g)$. Any  $g\in F_n$
represented by a loop meeting only strata below $H$ is elliptic. 

The
tree
$T$ is $\Phi
$-invariant, in the sense that $\ell(\Phi (\cc g))=\lambda \ell(\cc g)$, with  
$\lambda >1$ the Perron-Frobenius eigenvalue of $M_H$ as above.

As in [\LLD, \S3], we note:
\nom\polell
\thm{Lemma \sta} If $\cc g$ grows polynomially, then $g$ is elliptic in $T$
(it fixes a unique point).
\fthm

\demo{Proof}    $\lambda
^p\ell(\cc g) =
\ell(\Phi ^p(\cc g))$ is bounded by a constant times the word length of $\Phi
^p(\cc g)$, so $\ell(\cc g)=0$.
\cqfd\enddemo

Stabilizers of points of $T$ have rank $<n$, and  there are finitely many orbits of   points with non-trivial
stabilizer. These are general facts about trees with trivial arc stabilizers  [\GL]. In the case at hand, they may be deduced from the description of stabilizers used in Section 3.

\nom\refe
\thm{Definition \sta} If $\Phi $ is exponentially growing, let  $T$ be the $\Phi
$-invariant
\Rt{} constructed above. We let $m_i$,
$1\le i\le b$,  be representatives for the orbits of   points with non-trivial
stabilizer.  We denote by $G_i$ the
   stabilizer of $m_i$ (it is malnormal, but not always a free factor). It has
rank
$n_i<n$. The conjugacy classes of the
$G_i$'s are permuted by $\Phi $.
Replacing
$\Phi $ by a power, we may assume that $\Phi $ leaves $G_i$ invariant (up to
conjugacy) and therefore induces $\Phi _i\in\Out(G_i)$.
\fthm

\subhead Polynomial subgroups\endsubhead

Let $\Phi \in\Out(F_n)$. 
A subgroup $P$ is {\it polynomial\/}  if there exist $q\ge1$ and
$\alpha
$ representing $\Phi ^q$ such that  $\alpha (P)=P$ 
and $\alpha _{|P}$ is polynomially growing. If $P$ is polynomial, so are its
conjugates and its images by any automorphism representing a power of 
$\Phi
$. Furthermore, 
$\Phi
$ and all its   powers have the same polynomial subgroups (for negative
powers, recall that the inverse of a polynomially growing automorphism is
polynomially growing).
 
\nom\posub
\thm{Proposition \sta} Let $\Phi \in\Out(F_n)$.
\roster
\item The conjugacy class $\cc g$ of $ g\in F_n$ grows polynomially if and
only if
$g$ belongs to a polynomial subgroup.
\item Every non-trivial polynomial subgroup is contained in a
unique maximal one.
\item If $g\ne 1$ grows polynomially under $\alpha $
representing $\Phi $, the maximal polynomial subgroup $P(g)$ containing $g$ is the
set of elements growing polynomially under $\alpha $.
\item Maximal polynomial subgroups have finite rank, are malnormal, and
there are finitely many conjugacy classes of them.
\endroster
\fthm

Inequalities for the ranks of the polynomial subgroups will be given in Theorem
4.1. For the automorphism $\theta_n$ constructed in Section 5, the group $P_0$ is a maximal polynomial subgroup.

\demo{Proof} By induction on $n$. The results are true  if $\Phi $ is
polynomially growing, with $F_n$  the unique maximal polynomial
subgroup. If not, we consider $T$, $G_i$, $\Phi _i$ as in   \refe{}
(after replacing $\Phi $ by a power if needed).   Any  polynomial subgroup of
$G_i$ (relative to
$\Phi _i$) is a polynomial subgroup of $F_n$ (relative to $\Phi $).

Consider $g\ne 1$ such that  $\cc g$ grows polynomially (this is the case, in
particular, if
$g$ belongs to a polynomial subgroup). Let $P(g)$ be the set of $h\in F_n$
such that both $\cc h$ and $\cc{gh}$ grow polynomially. We prove by induction
on $n$ that
$P(g)$ is a polynomial subgroup.

Since $\cc g$ grows polynomially, the element $g$ fixes a point in $T$ by
Lemma
\polell{}, and by conjugating $g$ we may assume
$g\in G_i$. If
$h\in P(g)$, then $h$ is also elliptic. It fixes the same point as $g$ because
otherwise $gh$ would be hyperbolic in $T$, contradicting polynomial growth
of $\cc{gh}$. This shows
$P(g)\inc G_i$. By induction, $P(g)$ is a polynomial subgroup (relative to $\Phi
_i$, hence also to $\Phi $).

Clearly $P(g)$ contains $g$, as well as  every polynomial subgroup
containing
$g$: it is the largest polynomial subgroup containing $g$. This shows
assertions (1) and (2). 

Suppose $g$ and $\alpha $ are as in (3). As above, we may assume $g\in G_i$. 
Any element
growing polynomially under $\alpha $ is in $P(g)$. For the converse, note
that $\alpha (g)$ belongs to $P(g)$, hence to $G_i$. 
We deduce that  $G_i$ and $\alpha (G_i)$ have a
non-trivial intersection, hence are equal (they are point stabilizers in a tree
with trivial arc stabilizers).
By induction, any $h\in P(g)$ grows polynomially under $\alpha $. 

The finiteness statements in (4) are immediate by
induction. Let us prove malnormality.
Let $P$ be a maximal polynomial subgroup. If $P$ and $gPg\mi$  meet
non-trivially, they are equal (by the uniqueness statement in (2)), so we
only have to show that $P$ equals its normalizer. As above, we may assume
$P\inc G_i$. Since $P$ fixes a unique point $m_i$ in $T$, any   $g$ normalizing
$P$ belongs to $G_i$. By induction, $P$ equals its normalizer in $G_i$, so it
also equals its normalizer in $F_n$.
\cqfd\enddemo

\nom\ssfix
\thm{Lemma \sta} Given $\Phi \in\Out(F_n)$, the following are
equivalent:
\roster
\item There exists a polynomial subgroup of rank $\ge2$.
\item  There exist $q\ge1$ and $\alpha $ representing $\Phi ^q$ with
$\rk\fix\alpha \ge2$.
\endroster
\fthm

\demo{Proof} (2) obviously implies (1). For the converse, it suffices to
prove (2) under the assumption that
$\Phi $ is polynomially growing and $n\ge2$.  We distinguish two cases, as in
Definition \refer. First suppose
that some power of
$\Phi $ has a representative
$\alpha =\alpha _1*\alpha _2$. If $n_1$ or $n_2$ is $\ge 2$, we use induction.
Otherwise, $\alpha ^2$ is the identity of $F_2$ and (2) holds. In the second
case, we consider a decomposition $F_n=G_1*\langle t\rangle$. If $n\ge3$,
we use induction. If $n=2$, 
then $\alpha ^2\in\Aut(F_2)$ is of the form
  $a
\mapsto a, t\mapsto ta^i$ for some $i\in\Z$, and its fixed subgroup has
rank 2 (if $i\ne0$, it is generated by $a$ and $tat\mi$).
\cqfd\enddemo

\nom\deux
\thm{Corollary \sta} Let $\Phi \in\Out(F_n)$. Suppose that some non-trivial
conjugacy class $\cc g$ grows polynomially.  If $\cc g$ is
not  periodic, there exist
$q\ge1$ and
$\alpha $ representing $\Phi ^q$ with $\rk\fix\alpha  \ge2$.  If $\cc g$ is
   periodic, there exists $\alpha $ with $\rk\fix\alpha  \ge1$.
\fthm

\demo{Proof} Only the first assertion requires a proof. By Proposition \posub,
$g$ belongs to a polynomial subgroup.  It has rank at least 2, as otherwise
$\cc g$ would be periodic. We conclude by Lemma \ssfix.
\cqfd\enddemo

We also record the following easy fact:

\nom\fac
\thm{Lemma \sta} If every conjugacy class is $\Phi $-periodic, then $\Phi $ has
finite order. \cqfd
\fthm

In fact $\Phi $ has finite order as soon as   every conjugacy class of length
$\le2$ is periodic.

\head \sect The basic invariants \endhead

We write $x^+=\max(x,0)$.

  Given $\Phi  \in\Out (F_n)$, we now define   numbers $e,s,d,p,r,k$.
If we   consider several automorphisms simultaneously, we write
$e(\Phi ), s(\Phi )$,... so that no confusion arises. If we consider $\alpha
\in\Aut(F_n)$, we write $e,s,...$ for the invariants of the outer automorphism
$\Phi =\hat\alpha $  represented by $\alpha $.

$\bullet$ $e$ is the number of exponential strata of any improved
relative train track map representing a power of $\Phi $. It is also the
number of attracting laminations of $\Phi $ [\BFH, Subsection 3.1]. Note that $e$ is
an upper bound for the number of   growth types $(\lambda ,m)$
with $\lambda >1$ (see Theorem 6.2), and that $e=0$ if and only if $\Phi $ is
polynomially growing. Also note that $e(\Phi \mi)=e(\Phi )$ by Subsection 3.2 of
[\BFH].

$\bullet$ $s$ is the maximal length of a chain of attracting laminations
$\Lambda _0\supsetneq\dots\supsetneq
\Lambda _s$. Any   growth type $(\lambda ,m)$
with $\lambda >1$ satisfies $m\le s$ (see appendix).

$\bullet$ $d$ is the maximal degree of polynomial growth of conjugacy
classes. If $e=d=0$, then $\Phi $ has finite
order by Lemma \fac. One has $d(\Phi \mi)=d(\Phi )$ by [\Pig] (to be precise, one
also needs Proposition \posub{} to reduce to the polynomially growing case
studied in [\Pig], and Lemma 2.3 below to compare  $d$ to growth of elements).

$\bullet $ $p=(d-1)^+=\max(d-1,0)$.

$\bullet$ $r$ is the index of $\Phi $, computed using ranks of fixed
subgroups. Namely, we write $r_0(\alpha
)=(\rk\fix\alpha -1)^+=\max(\rk\fix\alpha -1,0)$ for $\alpha $ representing
$\Phi $. We then define
$r=\sum_i r_0(\alpha _i)$,  the sum being taken over a   set of
representatives of isogredience classes  (recall   that $\alpha ,\alpha '$
are isogredient if there is an inner automorphism conjugating them). 
One has
$r\le n-1$ by [\BH, Corollary 6.4] (see also [\GLL]). If $\rk\fix\alpha \ge2$, we say
that the isogredience class of $\alpha $ contributes $r_0(\alpha )$ to $r$. Only
finitely many classes contribute. 
 
$\bullet$ $k$ is the rank of the subgroup   generated by $\Phi $-periodic
conjugacy classes in the abelianization of $F_n$ (so $k\le n$). This number
plays an essential role in the proof of our main results (see the discussion in the
introduction and in Section 3).

The numbers  $e,s, d,p,k  $    do not change if $\Phi  $ is replaced by a
positive power (and $e,d,p, r,k$ do not change if $\Phi  $ is replaced by $\Phi \mi$; we do not know whether $s$ is the same for $\Phi  $ and $\Phi  \mi$).
The following lemma controls $r$, so that we can always    replace
$\Phi
$ by a power when  proving upper bounds for the invariants.

\nom\pui
\thm{Lemma \sta} $r(\Phi )\le r(\Phi ^q)$ for any $q\ge1$.
\fthm

\demo{Proof}   We have to
bound any finite sum $\sum_i r_0(\alpha _i)$ as above by $r(\Phi ^q)$. 
One always has $\rk\fix\alpha
\le\rk\fix
\alpha ^q$, since $\fix\alpha $ is a free factor of $\fix\alpha
^q$ by [\DS], so
 $r_0(\alpha _i)\le r_0(\alpha _i^q)$.
But this is not enough, because
$\alpha _i^q$ and
$\alpha _j^q$ may be isogredient for $i\ne j$ even though $\alpha _i,\alpha _j$
aren't.

We may assume that $\alpha _i^q$ and $\alpha _j^q$ are equal
if they are isogredient (by changing automorphisms within their isogredience
class), so we reduce to showing
$$\sum_{\alpha _i^q=\beta } r_0(\alpha _i)\le r_0(\beta )$$ for any $\beta $
representing $\Phi ^q$. We may also assume that $\fix\beta $ has rank
$\ge2$.

Note that $\fix \beta $ is $\alpha _i$-invariant and
contains $\fix\alpha _i$. We claim that the restrictions of the
$\alpha  _i$'s to $\fix\beta $ represent the same outer automorphism
$\tilde \Phi $ (of finite order), and define distinct   isogredience classes.
Assuming this, we conclude because
$$\sum_{\alpha _i^q=\beta } r_0(\alpha _i)\le r(\tilde \Phi )\le \rk\fix\beta
-1=r_0(\beta ).$$

To prove the claim, fix distinct $i,j$. We   write $\alpha _i=i_h\circ \alpha
_j$ (with
$i_h(g)=hgh\mi$) and we show $h\in\fix\beta $. We have
$\alpha _i^q=i_n\circ\alpha _j^q$ with $n=
h\alpha
_j(h)\dots\alpha _j^{q-1}(h)$. We deduce $n=1$, and this implies $\alpha
_j^q(h)=h$, so $h\in \fix\beta $ as required. This shows the first assertion. Next
(as in [\GLL, 2.4]), suppose that $i_m$ conjugates the restrictions of $\alpha _i$
and
$\alpha _j$, with
$m\in\fix\beta $. Then $\alpha _i\alpha _j\mi=i_{m\alpha _j(m\mi)}$ on
$\fix\beta $, hence everywhere because  $\alpha _i\alpha _j\mi$ is inner
and $\fix\beta $ has rank $\ge2$. This means that $i_m$ conjugates $\alpha
_i$ and $\alpha _j$ on the whole of $F_n$, so $\alpha _i,\alpha _j$
are isogredient, a contradiction to $i\ne j$.
\cqfd\enddemo

\nom\err
\example {Remark \sta}
By
Corollary \deux, $d(\Phi )>0$ implies that there exists $q\ge1$ with $r(\Phi
^q)>0$.
\endexample

The following lemma compares growth of conjugacy classes and growth of
elements.

\nom\trois
\thm{Lemma \sta} Let $\Phi \in\Out(F_n)$, and $d=d(\Phi )$.   Suppose that
some $g\in F_n$ grows polynomially with degree
$d'>d $ under iteration of a representative $\alpha
$ of $\Phi $. Then:
\roster\item
   $d'=d+1$.
\item
If    $ \fix\alpha  $ is non-trivial,  then $d'=1$ (and therefore $d=0$).
\endroster
\fthm

\example{Examples}  If $\alpha $ is a non-trivial inner automorphism, then
$d=0$ and any
 non-trivial $g$ grows linearly ($d'=1$).  In the   example $a\mapsto
bab\mi$,
$b\mapsto b^2ab\mi$, due to  Bridson-Groves, $\cc b$ grows linearly and $b$
grows quadratically (the fixed subgroup is trivial).\endexample

\demo{Proof}  We first prove (2).
Fix a basis of $F_n$, and $h\ne 1$ in $\fix\alpha $. Write
$\alpha ^p(g)=w_pg_pw_p^{\mi}$ with $g_p$ cyclically reduced.
The
length of  
$w_p$ grows with degree $d'$, because $g_p$ grows with degree at most
$d$.
We claim that the cancellation between the initial segments
of $w_p$ and $hw_p$ grows at most linearly. Assuming this for the
moment, we   consider
$\alpha ^p(hg)=hw_pg_pw_p^{\mi}$. If $d'>1$, the cyclic reduction of this
word grows with degree $d'$, so $\cc {hg}$ grows with degree
$d'$ under iteration of $\Phi $. This contradicts
$d'>d$.

To prove the claim, first note that the cancellation between $\alpha ^p(g)$
and $h^{\infty}=\lim_{p\to\infty}h^p$ grows at most linearly (see [\LLD, p\.
424]). The same holds for the cancellation between $w_p$ and $h^{\infty}$,
and therefore for that between $w_p$ and $hw_p$.

For (1), we let $\alpha _1$ be the restriction of $\alpha $ to  the polynomial 
subgroup
$P(g)$ (see item (3) in Proposition \posub). Let $\Phi _1$ be the outer
automorphism determined by $\alpha _1$. If $d(\Phi _1)=0$, then $\Phi _1$
has finite order by Lemma \fac, and $g$ grows at most linearly.  If not, we 
use   Corollary
\deux{}. After    replacing $\alpha $ by a power, we write
$\alpha _1=i_a\circ \beta $, where $a\in P(g)$ and $\beta \in\Aut(P(g))$ has
a non-trivial fixed subgroup.
Since $d(\Phi _1)>0$,    assertion (2) proved above implies that   
$a$ and
$g$ grow with degree at most
$d$ under iteration of $\beta $. Then   $$\alpha
^{n}(g)=[a\beta (a)\beta ^2(a)\dots\beta ^{n-1}(a)]\ \beta ^n(g)\ [a\beta
(a)\beta ^2(a)\dots\beta ^{n-1}(a)]\mi$$ grows with degree $\le d+1$.
\cqfd\enddemo

Our next goal will be to understand the behavior of the invariants in the
induction process, so we consider $\Phi $ as in \refer{} or \refe.

First suppose $\alpha =\alpha _1*\alpha _2$ as in the first item of
Definition \refer. 
We denote by
$k_i,r_i,d_i$ the invariants associated to
$\hat\alpha _i$.

\nom\kb
\thm{Lemma \sta}
\roster
\item If $\fix\alpha _1$ or $\fix\alpha _2$ is trivial, then $r=r_1+r_2$ and
$\max(d_1,d_2)\le d\le 1+\max(d_1,d_2)$.
\item If $\fix\alpha _1$ and $\fix\alpha _2$ are both non-trivial, then
$r=r_1+r_2+1$ and $d$ equals 1 or $\max(d_1,d_2)$.
\item $k=k_1+k_2$.
\endroster
\fthm

\demo{Proof}  The
assertions about $r$ are proved in [\GLL]. Here is the idea. If an isogredience
class other than that of $\alpha $  contributes to $r$, it contributes the same
amount to either
$r_1$ or
$r_2$ (but not both). The class of $\alpha $ contributes $(\rk\fix\alpha _i-1)^+$
to $r_i$ and $(\rk\fix\alpha  -1)^+$ to $r$. To conclude, note that $\fix\alpha
=\fix\alpha _1*\fix\alpha _2$.

 The number $d$ is bounded from
below by $\max(d_1,d_2)$, and from above by the maximal degree of growth of
elements of
$G_i$ under $\alpha _i$, so the assertions about $d$ follow from Lemma
\trois. The assertion about
$k$ is easy and left to the reader.
\cqfd\enddemo

We now consider a
decomposition
$F_n=G_1*\langle t\rangle$  and a representative $\alpha $ of $\Phi $ such
that $\alpha (G_1)=G_1$ and $\alpha (t)=tu$ with $u\in G_1$, as in item (2) of
Definition \refer. Let
$k_1, r_1,d_1$ be the invariants associated to the outer automorphism $\Phi
_1$ of
$G_1$ represented by
$\alpha _1=\alpha _{|G_1}$.

\nom\ka
\thm{Lemma \sta}  \roster
\item
$r\le r_1+1$, with equality if and only if both $\alpha
_1$ and
$i_u\circ\alpha _1$ have non-trivial fixed subgroups.
\item
$ d\le
d_1+1$.
\item   $k\le k_1+1$, with equality  if and only there
exists
$a\in G_1$ such that $ta$ is $\alpha $-periodic.

\endroster
\fthm

\demo{Proof} The assertion about $r$ is proved in [\GLL].  One now has $\fix\alpha
=\fix\alpha _1*t\fix(i_u\circ\alpha _1)t\mi$. The contribution of the
isogredience class of
$\alpha
$ to
$r$ is  
$(\rk\fix\alpha _1+\rk\fix (i_u\circ\alpha _1)-1)^+$, whereas its contribution to
$r_1$ is 
$(\rk\fix\alpha _1 -1)^+ +(\rk\fix
(i_u\circ\alpha _1)-1)^+$.

Let us prove (2). As in
the proof of Lemma \trois, replace $\alpha $ by a power and write $\alpha
_1=i_a\circ
\beta _1
$, where
$a\in G_1$ and
$\fix\beta _1$ is non-trivial, using Corollary \deux; if $d_1=0$, we take
$\beta _1$ to be the identity, using Lemma \fac.  By
Lemma \trois, all elements of $G_1$ grow with degree at most $d_1$ under
$\beta _1$.
Now 
$\beta =(i_a)^{\mi}\circ\alpha $ sends $t$ to an element of the form  $u'tu''$,
so $t$ grows with degree at most $d_1+1$ under $\beta $. This implies $d\le
d_1+1$.

We now study
$k$.
Represent conjugacy classes of
$F_n$ as cyclic words whose letters are either $t^{\pm1}$ or non-trivial
elements of
$G_1$. Let
$\ov w$ be a
cyclic word not contained in $G_1$.  Split it before every letter $t$ and after 
every
$t^{\mi}$. This expresses $\ov w$  as a  product of subwords of the
form
  $tat^{\mi}$,
$ta$,
$at^{\mi}$, $a$, with $a\in G_1$ and $a$ non-trivial in subwords
$tat^{\mi}$ or $a$.
 This
decomposition  of $\ov w$ is preserved   by
$\alpha
$, so that
$\ov w$ is $\Phi $-periodic if and only if every subword is an $\alpha
$-periodic element.

   Let  $K$ be the quotient of the abelianization of $F_n$ by the
subgroup generated by all $\Phi _1$-periodic conjugacy classes (of elements
of $G_1$). It has torsion-free rank $n-k_1$. To
bound $k$, we need to control the image  in $K$  of $\alpha $-periodic
elements of the form
$tat^{\mi}$ or  $ta$ (we treat $at^{\mi}$ as the inverse of 
$ta^{\mi}$).

If
$tat^{\mi}$ is $\alpha $-periodic, then $a$ is $(i_u\mi\circ \alpha _1)
$-periodic, so $\ov a$ is
$\Phi _1$-periodic  and 
$tat^{\mi}$ maps trivially to $K$. On the other hand, the image of $ta$ in $K$
has infinite order. But, if
$ta$ and $ta'$ are both $\alpha $-periodic, then $a^{\mi}a'$ is $\alpha
_1$-periodic and therefore
$ta$,
$ta'$ have the same image in $K$. This proves the lemma.
\cqfd\enddemo

We now suppose that $\Phi $ has exponential growth, and we use the \Rt{} $T$
as in \refe.  
Recall that we have defined  invariant subgroups $G_i$ and 
restrictions   $\Phi _i$.
The following lemma allows us to bound the invariants of $\Phi $ in
terms of those of the $\Phi _i$'s (which we denote with   the subscript
$i$). All sums are over $i$.
\nom \ine
\thm{Lemma \sta} The invariants of $\Phi $ satisfy:
$$  \gather e\le 1+\sum e_i\\
 s\le 1+\max s_i\\
d=\max d_i\\
r=\sum  r_i\\
k\le\sum k_i.
\endgather $$
\fthm

\demo{Proof} The results for $d$ and $k$ are direct consequences of
Lemma
\polell{} (note that for $g\in G_i$   the
growth of $\cc g$ is the same  for
$\Phi _i$ as for $\Phi $). The equality for $r$ is proved in [\GLL]: if an isogredience
class contributes to $r$, it contributes the same amount to exactly one $r_i$. 

As in the
previous section,    let
$H$ be the highest exponential stratum of the train track, and let $H'$  be
the union of
$H$ with all the strata above it. Every exponential stratum $\tilde H$ other
than
$H$ is contained in a component $\tilde G$ of the closure of $G\setminus H'$. It
follows from the way $T$ was constructed (see [\GJLL] and the previous
section) that 
the fundamental group of
$\tilde G$ fixes a point   of $T$, whose stabilizer is conjugate to some $G_i$.
There is a bijection between exponential strata and attracting laminations
[\BFH, Definition 3.1.12]. By [\BFH, Definition  3.1.5], the lamination associated to
$\tilde H$ is an attracting lamination of
$\Phi _i$. This shows the bound for $e$. If  $\Lambda
_0\supsetneq\dots\supsetneq
\Lambda _s$
 is a chain of laminations, then $\Lambda _1,\dots,
\Lambda _{s }$ are attracting laminations of some $\Phi _i$, so $s-1\le s_i$.
\cqfd\enddemo

\head \sect Bounding ranks of stabilizers \endhead

In this section we assume that $\Phi $ is not polynomially growing, and we
consider the invariant
\Rt{}
$T$ as in \refe.

\nom\inee
\thm{Proposition \sta} Let $n_i$ be the rank of the stabilizer $G_i$, and $k_i=k(\Phi _i)$. Then: 
$$  \gather
\sum_{i=1}^b(n_i-1)\le n-2 \tag1\\
\sum_{i=1}^b(3n_i-2-k_i) \le 3n-6-k  .\tag2
\endgather $$
 If equality holds in (1), there is only one exponential stratum.
 \fthm

Recall that $k\le \sum k_i$, so (2) is implied by the simpler inequality $\sum_i(3n_i-2)\le 3n-6$.  But this simpler inequality is not always true, for instance if $\Phi$ is induced by a pseudo-Anosov homeomorphism of a punctured torus or a four-punctured sphere. This explains the introduction of   the invariant $k$.

\demo{Proof of Proposition \inee}

Recall that $T$ was constructed using the highest exponential stratum $H$ of $f:G\to G$.
A description of the $G_i$'s is provided by Theorem 6.0.1  of [\BFH], in terms of a subgraph $Z\inc G$. 

$\bullet$ First suppose that $H$ is the highest stratum. Then $Z$ is the union of all strata below $H$. 

If $H$ is not geometric (in the sense of Definition 5.1.4 of [\BFH]), then the $G_i$'s are (up to conjugacy) the fundamental groups of the non-contractible components of $Z$. In particular they are free factors and $\sum n_i\le n$. If $b=0$, we have $n\ge3$ since every automorphism of $F_2$ is induced by a homeomorphism of a punctured torus. 
If
$b=1$ or
$b=2$, the existence of the exponential stratum $H$ prevents 
$n_1,n_2$ from being too big. More precisely, we have   $n_1< n-1$ if $b=1$
by [\BFH, Lemma 3.2.1], and similarly $n_1+n_2<n$ if $b=2$.

We get $\sum (n_i-1)\le n-3$ and $\sum (3n_i-2) \le 3n-6$, with equality in the second inequality possible only if $b=0$ or $b=3$. The proposition is true  in this case since $\sum k_i\le k$.

If $H$ is geometric, the components of $Z$ account for all $G_i$'s but one. The exceptional one, say $G_b$, is cyclic, generated by the homotopy class of a loop $\rho$ (an indivisible Nielsen path) based at a point $v\notin Z$.  If $b $ is different from 1 and 4, the previous argument yields $\sum_{i<b} (n_i-1)\le n-3$ and $\sum _{i<b}(3n_i-2) \le 3n-7$, and the proposition holds since $n_b=1$. 

If $b=1$, we have $\sum(n_i-1)=0\le n-2$. On the other hand, $\sum (3n_i-2)=1$. It  is bounded by $3n-6$, except if $n=2$. But in this case  $\Phi$ is induced by  a homeomorphism of a punctured torus, so $k=0$ (whereas $k_1=1$). We get    $3n_1-2-k_1=0=3n-6-k$.

When $b=4$, we have $\sum (n_i-1)\le n-3$ and $\sum _{i\le 4}  (3n_i-2)\le3\sum _{i<4} n_i-5$.
If $\sum _{i<4} n_i<n$, we are done. If not, we prove (2) by showing $k<\sum _{i\le 4} k_i$. 

The geometric stratum $H$ is associated to  a connected surface $S$, as in Definition 5.1.4 of [\BFH]. This surface has four boundary components $C_i$, with  $C_4$   identified to $\rho$. 
It  is a punctured sphere because $\sum _{i<4} n_i=n$. The  $C_i$'s  represent $\Phi$-periodic conjugacy classes $z_i$, and the relation $\sum z_i=0$ holds in the abelianization of $F_n$. Since $z_4$ may be expressed in terms of $z_1,z_2,z_3$, we get  $k\le\sum_{i<4} k_i=(\sum_{i\le 4} k_i )-1$.

Note that $\sum (n_i-1)= n-2$ is possible only if $n=2$ (and then the exponential stratum is unique).

$\bullet$ Now suppose that there are non-exponential strata above $H$. Then $Z$ contains all strata below $H$, no edge of $H$, and possibly edges from strata above $H$.

Consider the union of $H$ and all strata below it. By aperiodicity, this subgraph has a component containing $H$, we call it $Y^0$. Then define connected subgraphs $Y^0\inc Y^1\inc\dots \inc Y^q=G$ such that $Y^j\setminus Y^{j-1}$ contains exactly one edge. 

If $H$ is not geometric,  the $n_i$'s are the Betti numbers of the non-contractible components of $Z$. 
Define $Z^j=Z\cap Y^j$. Let $n^j_i$ be the Betti numbers of the non-contractible components of $Z^j$, and $n^j$ the Betti number of $Y^j$. Since $Y^0$ is $f$-invariant, the previous argument yields $\sum (n^0_i-1)\le n^0-3$ and $\sum (3n^0_i-2) \le 3n^0-6$.   Using induction on $j$, one shows    $\sum (n^q_i-1)\le n^q-3$ and $\sum (3n^q_i-2) \le 3n^q-6$ (when passing from $j$ to $j+1$, the left hand sides cannot increase more than the right hand sides).   This proves the proposition since $Z^q=Z$.

If $H$ is geometric, we again have to consider $\rho$ and $v$. We define $Z^j=(Z\cap Y^j)\cup\{v\}$. The point $v$ is an isolated point of $Z^0$. The numbers $n^j_i$ and $n^j$ are defined as before, except that we add 1 to the Betti number of the component  $Z^j_v$ of $Z^j$ which contains $v$ (in particular, we always consider it as non-contractible). The ranks $n_i$ of the groups $G_i$ are the $n^q_i$'s (the exceptional group $G_b$ is generated by $\pi_1(Z^q_v)$ and the class of $\rho$).

The inequalities $\sum (n^j_i-1)\le n^j-2$ and $\sum (3n^j_i-2) \le 3n^j-5$ are true for $j=0$, hence for $j=q$ by induction. If $\sum (3n^j_i-2) = 3n^j-5$  holds for $j=q$, it holds for $j=0$ so there is  a 
punctured sphere $S$ as above. 
Let 
$z$ be the conjugacy class   represented by $\rho$.  Removing  $z$ from the set of periodic conjugacy classes does not change the subgroup generated in the abelianization of $F_n$. In other words, $z$ does not contribute to $k$. But it contributes to $k_b$ since it generates a free factor in $G_b$.
We get $k<\sum k_i$, and the proposition is proved. 
\cqfd\enddemo

\head \sect Proof of the main results\endhead

In this section we give upper bounds for the invariants $e,s,d,p,r$ introduced
in Section 2.  Recall that we are free to replace $\Phi $ by a power  (see
Lemma \pui{} and the paragraph preceding it).

\nom\pol
\thm{Theorem \sta} Given $\Phi \in\Out(F_n)$, let $P_j$ be
representatives for the  conjugacy classes of maximal
polynomial subgroups. Then $$\aligned e+\sum_j (\rk P_j-1)^+&\le n-1\\
4e+k+2\sum_j (\rk P_j-1)^+&\le 3n-2.\endaligned$$
\fthm

Recall that $x^+=\max(x,0)$.

\demo{Proof} The result is trivially true if $\Phi $ is polynomially growing
(with
$e=0$ and $k\le n$), so we replace $\Phi $ by  a power and we consider $T$,
$G_i$,
$n_i$,
$\Phi _i$ as in
\refe. 

If the action on $T$ is free, then by Lemma \polell{} there is no non-trivial
polynomially growing conjugacy class (so $k=0$) and  no non-trivial polynomial
subgroup. Furthermore, 
$e=1$. The theorem is  true in this case since $n\ge2$, so we assume $b>0$. We
argue by induction on $n$.

We may assume that each $P_j$ is
contained in a (unique) $G_{i_j}$. The groups $P_j$ contained in a given
$G_i$ are non-conjugate maximal polynomial subgroups of $\Phi _i$. 

 Successively using Lemma \ine{}, the induction hypothesis, and  Proposition
\inee{}, we now write:
   $$\aligned e+\sum_j (\rk P_j-1)^+&\le 1+\sum_{i=1}^b
e_i+\sum_{i=1}^b\ \sum_{i_j=i}(\rk P_j-1)^+\\
&\le 1+\sum_{i=1}^b(n_i-1)\\
&\le n-1 \endaligned$$
and
$$\aligned 4e+k+2\sum_j (\rk P_j-1)^+&\le 4+4\sum_{i=1}^b
e_i+k+2\sum_{i=1}^b\ \sum_{i_j=i}(\rk P_j-1)^+\\
&\le 4+k+\sum_{i=1}^b(3n_i-2-k_i)\\
&\le 4+3n-6\\  &\le
3n-2.\endaligned$$
\cqfd\enddemo

\nom\rpol
\thm{Proposition \sta} If $\Phi \in\Out(F_n)$ is polynomially growing, then
$p+r\le n-1$.
\fthm

\demo{Proof} By induction on $n$. After replacing $\Phi $ by a power, we  may
assume that we are in the situation of Definition
\refer.
 There are two possibilities.

$\bullet$ First suppose that  some $\alpha $
representing $\Phi $ preserves a  non-trivial decomposition
$F_n=G_1*G_2$. Let $n_i$ be the rank of $G_i$. Define $\alpha _i,\Phi _i,
d_i, p_i, r_i$ by considering the automorphism induced on $G_i$.
Recall from Lemma \kb{} that
$r=r_1+r_2+1$ if both $\fix\alpha _1$ and $\fix\alpha _2$ are non-trivial,
$r=r_1+r_2$ otherwise. Furthermore,   $d\le\max(d_1,d_2)+1$
and $p\le\max(p_1,p_2)+1$.

If $r=r_1+r_2$, we write $$p+r\le p_1+p_2+1+r_1+r_2\le
n_1-1+n_2-1+1=n-1,$$
so we assume $r=r_1+r_2+1$. 
By Lemma \kb, we have $d=1$ or
$d=\max(d_1,d_2)$. In both cases $p=\max(p_1,p_2)$ and
$$p+r\le p_1+p_2 +r_1+r_2+1\le
  n-1.$$

$\bullet$ Now suppose that there is a
decomposition $F_{n}=G_1*\langle t\rangle$ and a representative
$\alpha
$ of $\Phi $ which leaves $ G_1\simeq F_{n-1} $ invariant and maps   $t$ to
$tu$ with
$u\in G_1$. As in \refer, we   assume that $u$ cannot be written as $a\alpha
(a\mi)$ with $a\in G_1$. In particular, $u\ne 1$. Define
$\alpha _1,\Phi _1, d_1, p_1,r_1$ by restricting to $G_1$. Recall (Lemma \ka)
that
$r\le r_1+1$, with equality if and only if both $\alpha _1$ and
$i_u\circ\alpha _1$ have non-trivial fixed subgroups, and $p\le
p_1+1$.

We show that $r=r_1+1$ implies $p=p_1$, so that in all cases $$p+r\le
p_1+r_1+1\le n_1-1+1=n-1.$$
Suppose $r=r_1+1$. Since $\fix(i_u\circ\alpha _1)$
is non-trivial, there is a non-trivial $\alpha $-fixed element of the form
$txt\mi$, and $u$ is $\alpha $-fixed by property ne-(iii) of [\BFH, Theorem
5.1.5]. In particular, the element $t$ grows linearly under $\alpha $. If $d\ge2$,
some element or conjugacy class of $G_1$ grows with degree   $d$, and
$d=d_1$  by Lemma \trois. If $d\le1$, we have $p=p_1$.
\cqfd\enddemo

Combining Theorem \pol{} and Proposition   \rpol, we deduce:

\nom\pine
\thm{Theorem \sta} Given any $\Phi \in\Out(F_n)$, we have:
$$\aligned e+p+r&\le n-1\\
4e+2p+2r+k&\le 3n-2.\endaligned$$
\fthm

\demo{Proof} By Theorem \pol, it suffices to show $p+r\le \sum_j (\rk P_j-1)^+$.
We may replace $\Phi $ by a power, and so assume that each $P_j$ is fixed (up to
conjugacy). Since $P_j$ equals its normalizer, there is a well-defined induced
$\Phi _j$, with   associated $p_j$ and $r_j$, and $p=\max p_j$.

Let $\alpha $ be a representative of $\Phi $. 
If $\fix\alpha $ is non-trivial, it is contained in
a conjugate of some $P_j$, and by changing $\alpha $ in its isogredience
class we may assume $\fix\alpha \inc P_j$. This implies that $P_j$ is
$\alpha $-invariant (by uniqueness of maximal polynomial subgroups).
Furthermore, if
$\fix\alpha
$ and
$\fix\alpha '$ are contained in $P_j$, 
then the  restrictions of $\alpha $ and $\alpha '$ to $P_j$ represent the same
outer automorphism (if $\alpha '=i_h\circ\alpha $, then $P_j$ is $i_h$-invariant,
so
$h\in P_j$
because $P_j$ equals its normalizer); if $P_j$ has rank
$\ge2$, and $\alpha ,\alpha '$ are not
isogredient, the restrictions    are not isogredient (see [\GLL] or the proof of
Lemma \pui).

Write  $r=\sum  r_0(\alpha _m)$, the sum being taken over non-isogredient
automorphisms
$\alpha _m$. We may assume
$\fix\alpha _m\inc P_{j_m}$, with $P_{j_m}$ of rank $\ge2$. Now $$  p+r\le\sum
_j p_j+\sum_j\sum_{\fix\alpha _m\inc P_j} r_0(\alpha _m)\le\sum
_j p_j+\sum_j r_j\le \sum (\rk
P_j-1)^+$$
by Proposition \rpol.
\cqfd\enddemo

The inequality   $4e+2p+2r+k\le 3n-2$ may be an equality,   for instance
for the automorphism $\alpha $ of $F_3$ defined by
$a\mapsto a$,
$b\mapsto b$,
$c\mapsto aca\mi$, with $n=k=3$, $d=1$, $r=2$. Note that $r(\hat\alpha )=2$,
but
 the contribution of $\alpha $ is $r_0(\alpha )=(\rk\fix\alpha -1)^+= 1$, as
$\fix\alpha $ has rank 2. The other non-zero contribution is from   $\alpha
'=i_{a\mi}\circ\alpha $, whose fixed subgroup also has rank 2. The point of the
next theorem is that this is a general phenomenon:   if $\Phi \in\Out(F_n)$
satisfies $4e+2p+2r+k= 3n-2$, and
$d>0$, then
$r$   has to be carried by at least two isogredience classes (no
representative $\alpha $ of $\Phi $ satisfies $\rk\fix\alpha =r+1$).

\nom\inal
\thm{Theorem \sta} Given any $\alpha \in\Aut(F_n)$, we have
$$4e+2d+2\rk\fix\alpha +k\le 3n+1.$$
\fthm

Note that Theorem \pine{}, together with the inequalities $\rk\fix\alpha \le
r+1$ and
$d\le p+1$, implies $4e+2d+2\rk\fix\alpha +k\le 3n+2$.

\demo{Proof} It suffices to show
$$4e+2p+2\rk\fix\alpha +k= 3n\implies d=0,\tag*$$
 since any $\alpha $ satisfying $4e+2d+2\rk\fix\alpha +k= 3n+2$ must satisfy
$d=p+1$ and therefore
$4e+2p+2\rk\fix\alpha +k= 3n$.  The proof of $(*)$ is by induction on $n$.

Consider  $\alpha $ satisfying  $4e+2p+2\rk\fix\alpha +k= 3n$. Note that  
$\fix\alpha $ must have   rank $r+1$. 
We argue by way of contradiction, assuming $d>0$. By Remark \err, we may
also
assume
$r>0$ (we are free to replace $\alpha $ by a power since $e,p,k$ do not change
and
$\rk\fix\alpha
$ may only increase). In particular, $\fix\alpha $ has rank $r+1\ge2$.

First suppose that $\alpha $ is polynomially growing. Note that in this case
$4e+2p+2\rk\fix\alpha +k= 3n$ together with  $ d=0$  imply that $\alpha $ is  
the identity, since $\alpha $ has a power which is inner by Lemma \fac{}, and
$\rk\fix\alpha =k=n$.

After raising $\alpha $ to a power, we are in the situation of    
  Definition \refer.  The representative of $\hat\alpha $ introduced in
Definition
\refer{} is not necessarily equal to $\alpha $, so we denote it by $\beta $.   As
usual, we distinguish two cases.

$\bullet$ First suppose that  there is a decomposition
$F_n=G_1*G_2$ such that   $\beta $   leaves
each $G_i$ invariant.  The main difficulty is to show that we may assume
$\beta =\alpha $.

  Consider the action of
$F_n$ on the simplicial tree $T_0$ associated to the free product $G_1*G_2$.
This action is $\hat\alpha 
$-invariant, and there is an isometry $H$ of $T_0$ representing $\alpha $ in
the sense that $\alpha  (g)H=Hg$ for all $g\in F_n$ (see [\GLL]). We have to
show that
$H$ fixes an edge $e$: if it does, the stabilizers $G'_1,G'_2$ of the endpoints of
$e$ are
$\alpha $-invariant, and $\alpha $ preserves the free product decomposition
$F_n=G'_1*G'_2$.

Since $\fix\alpha  $ has
rank $\ge2$, the map $H$ has a fixed point [\GLL, 1.1]. The stabilizer of this
point is $\alpha $-invariant, so  we may assume $\alpha (G_1)=G_1$
(possibly exchanging the roles of $G_1$ and $G_2$ and replacing $\alpha $
by an isogredient automorphism). By an argument given in the proof of \pine,
$\alpha $ and
$\beta
$ induce the same outer automorphism $\Phi _1$ of $G_1$ because $G_1$
equals its normalizer. The numbers $k_i, p_i, r_i$ used below refer to $\Phi
_1$ and to 
$\Phi _2=\hat{\beta _{|G_2}}$

First suppose $\fix\alpha \inc G_1$ (this holds in particular if $H$ has a
unique fixed point). Then 
$$\align
3n&=  2\rk\fix\alpha +k+2p\\&\le  2\rk\fix\alpha
_{|G_1}+k_1+k_2+2p_1+2p_2+2
\\ &\le3n_1+3n_2-2+2\\&=3n
, \endalign  $$
using Lemma \kb{} and Theorem \pine.

All inequalities must be equalities. In particular $d(\Phi _1)=0$ because $(*)$
is true on $G_1$ by induction, and as pointed out earlier this implies that 
 $\alpha $ is the identity on
$G_1$.
  Also $2p=2p_1+2p_2+2$, so $d\ge2$. But $r_2=0$, because
$r_1+r_2\le r=\rk\fix\alpha
-1=r_1$ since $\fix\alpha =  G_1$.  
  Remark \err{} implies $d_2=0$, contradicting Lemma \kb.

This proves that $\fix\alpha \inc G_1$ cannot hold. Therefore $H$ fixes an
edge, so we may indeed assume
$\alpha =\alpha _1*\alpha _2$.  Then $\fix\alpha =\fix\alpha _1*\fix\alpha
_2$. Since we have ruled out $\fix\alpha \inc G_i$,
 both $\fix\alpha _1$ and
$\fix\alpha _2$ are nontrivial.
As in the proof of Theorem
\rpol, this implies
$p=\max(p_1,p_2)$. 

We now write:
$$\align
3n&=  2\rk\fix\alpha +k+2p\\
&\le 2\rk \fix\alpha _1+2\rk \fix\alpha _2
+k_1+k_2+ 2p_1+2p_2
\\ &\le3n_1+3n_2\\&\le3n
. \endalign  $$

Applying $(*)$ to $\alpha _1$ and $\alpha _2$ (using the induction
hypothesis), we find that $\alpha _1$ and $\alpha _2$ are the identity, so
$\alpha
$ is the identity.

$\bullet$ Now consider a
decomposition $F_n=G_1*\langle t\rangle$ as in Definition \refer, with
$\beta $ representing $\hat\alpha $ such that 
$\beta (G_1)=G_1$ and $\beta (t)=tu$. The numbers $k_1,p_1$ refer
to
$\hat{\beta _{|G_1}}$. We have shown the inequality
$p+r\le p_1+r_1+1$ in the proof of Proposition
\rpol. Using Theorem \pine, we deduce
$$\align3n-2&=  2r+k+2p\\
&\le  2r_1 +k_1+2p_1+2+k-k_1\\
&\le
3n_1-2+2+k-k_1\\
&\le3n-3+k-k_1,\endalign$$
showing $k\ge k_1+1$. By Lemma \ka, some decomposition
$F_n=G_1*\langle ta\rangle$ is invariant under a power of $\beta $, and we
reduce to the  previous case.

The   proof of Theorem \inal{} is now complete for polynomially growing
automorphisms. In the general case, we consider the invariant {}\Rt{}
$T$ as in \refe.

All elements growing polynomially under $\alpha $
are contained in an
$\alpha $-invariant  stabilizer, say $G_1$. Let $\alpha _1=\alpha _{|G_1}$,
so that $\fix\alpha =\fix\alpha _1$ and $r(\hat\alpha )=r(\hat\alpha _1)$. All  
$\Phi _i$'s with
$i\ge2$ have
$r_i=0$ by Lemma \ine, hence
$d_i=0$ by Remark \err. Since $d=\max d_i$, it suffices to show $d_1=0$.

Writing
$$\align 4e +2\rk\fix\alpha +k+2p&\le4+4\sum e_i+2\rk\fix\alpha
_1+k-\sum k_i +\sum k_i+2p_1\\
&\le (4e_1+2\rk\fix\alpha _1+k_1+2p_1)
+\sum_{i>1}(4e_i +k_i)+4+k-\sum_{i\ge1} k_i\\
&\le (4e_1+2\rk\fix\alpha _1+k_1+2p_1)
+\sum_{i>1}(3n_i-2)+4+k-\sum_{i\ge1} k_i
\\&\le (4e_1+2\rk\fix\alpha _1+k_1+2p_1)+3n-3n_1,\endalign$$
we get $4e_1+2\rk\fix\alpha _1+k_1+2p_1=3n_1$, and $d_1=0$ by
the induction hypothesis.
\cqfd\enddemo

By Corollary \deux, any $\Phi $ with $d>0$ has a power $\Phi ^q$ represented by an
automorphism $\alpha $ with $\rk\fix\alpha \ge2$. In particular, $\Phi ^q$
satisfies 
$d\le p+r$. We thus get from Theorems \pine{} and \inal:

\nom\dcor
\thm{Corollary \sta} Given $\Phi \in\Out(F_n)$, we have
$$\aligned
e+d&\le n-1\\
4e+2d&\le 3n-2\\
4e+2d&\le 3n-3\quad\text{if $d>0$}. 
\endaligned$$ \cqfd
\fthm

As is easily checked, this corollary is   equivalent to saying that
$(e,d)$ belongs to the closed quadrilateral with vertices 
$(0,0),(0,n-1),(\frac{n-1}2,\frac{n-1}2), (\frac{3n-2}4,0)$ pictured on Figure 1.

We also get:

\nom\pcor
\thm{Corollary \sta} Given $\alpha \in\Aut(F_n)$, we have
$$\aligned
e+(d-1)^++\rk\fix\alpha &\le n  \\
4e+2d+2\rk\fix\alpha &\le 3n+1\qquad \text{($\le 3n$ if $d=0$)}.
\endaligned$$ \cqfd
\fthm

We finally prove:

\nom\cha
\thm{Theorem \sta}  If $\Phi \in\Out(F_n)$ is not polynomially growing, one has
$2s+p+r\le n-2$.   In particular, $\ds s\le\frac
n2-1$.
\fthm

It follows from the appendix that $m\le s$ if a conjugacy class grows like
$\lambda ^pp^m$, with $\lambda >1$, under iteration of $\Phi $.  The optimality of 
$ s\le\frac
n2-1$ will be shown in the next section.

\demo{Proof}
The result is true if $s=0$, by Theorem
\pine, so we assume $s>0$. 
We argue
by induction on $n$, using   the invariant \Rt{}    $T$ as in \refe.  Since
$s>0$,  there is at least one
$i_0\in\{1,\dots,b\}$ such that $\Phi _{i_0}$ is not polynomially growing. By 
Lemma
\ine,  there is such an ${i_0}$  with
$s\le 1+s_{i_0}$.  

We then
get 
$$ 2s+p+r\le2s_{i_0}+2+\sum _{i=1 }^bp_i+\sum _{i=1 }^br_i\le
n_{i_0}+\sum_{{i\ne {i_0}}}(n_i-1)=1+\sum_{i=1 }^b(n_i-1)
$$
by the induction hypothesis and Theorem \pine. Since $s>0$, Proposition \inee{} yields  
 $\sum_{i=1 }^b(n_i-1)<
n-2$ and the theorem is proved.
\cqfd\enddemo

\thm{Corollary \sta} If $\alpha \in\Aut(F_n)$ is not polynomially growing,  one
has
$2s+d+\rk\fix\alpha \le n$ and $2s+d\le n-2$. 
\fthm

\demo{Proof}
The first inequality is clear. The second follows from Remark \err.
\cqfd\enddemo

\head\sect Examples\endhead

We give examples, and we show that the inequalities of Corollaries \dcor{} and
\pcor{} are optimal.

\subhead Automorphisms of $F_2$\endsubhead

Any $\Phi \in\Out(F_2)$ is induced by a homeomorphism of a punctured
torus. Some power of $\Phi $ is either the identity, or a Dehn twist, or a 
  pseudo-Anosov map.

The simplest pseudo-Anosov automorphism is $a\mapsto ab$, $b\mapsto a$.
For future reference,   we note that its square $\tau $, which sends   $a$ to
$aba$ and
$b$ to $ba$, fixes  the
commutator
$[a,b]=aba\mi b\mi$.

\subhead  Geometric automorphisms
\endsubhead

An automorphism $\Phi \in\Out(F_n)$ is geometric if it is induced by a
homeomorphism of a compact surface $\Sigma $ with fundamental group $F_n$.
For
$\Phi $ geometric, it follows from Nielsen-Thurston theory that the growth
  of any non-periodic conjugacy class under iteration of $\Phi  $ (always 
equivalent to some $\lambda ^pp^m$)
 is either linear 
or purely exponential:
$(\lambda ,m)=(1,1)$, or $\lambda >1$ and
$m=0$. 

We now construct a geometric automorphism $\Phi _n$ of $F_n$ with $e$
equal to the maximal value $e_n=\left[\frac{3n-2}4\right]$. It has a 
representative $\varphi _n\in\Aut(F_n)$ with non-trivial fixed subgroup. 

 Write $n=4\ell+3+\delta $ with $0\le\delta \le3$ (for $n=2$, we take $\Phi
_2=\hat\tau $). Construct a compact surface
$\Sigma _n$ with fundamental group $F_n$ by gluing $2\ell +\delta $
once-punctured tori and $\ell+1$ four-punctured spheres (see Figure 2 for a
picture with
$\delta =3$). The number of subsurfaces is $e_n=3\ell +\delta
+1=\left[\frac{3n-2}4\right]$.
  Consider an orientation-preserving
homeomorphism of $\Sigma _n
$ inducing a pseudo-Anosov map on each of the subsurfaces. The induced
automorphism $\Phi _n
\in\Out(F_n)$ satisfies 
$e(\Phi _n)=e_n$. 

If all pseudo-Anosov maps used in the
construction have distinct dilation factors $\lambda _i$, then $\Phi _n$ has
$e_n$ different exponential growth types $(\lambda _i,0)$.

This example explains the appearance of  
four-punctured spheres in the proof of Proposition \inee.

\example{Remark \sta} If $\Phi \in\Out(F_n)$ is induced by a homeomorphism of 
a compact orientable surface of genus $g$ with $b$ boundary components, one has
$n=2g+b-1$ and
$k\ge b-1$, so Theorem
\pol{} yields $e\le \frac{3g+b-2}2$. It is easy to see that this bound is optimal.
\endexample

\subhead Nested laminations\endsubhead

The inequality $\ds s\le\frac
n2-1$ of Theorem \cha{} is an equality   for the automorphism 
 $\alpha $ of $F_{2\ell}$ defined  by:

$$\left\{\aligned
a_1&\mapsto a_1b_1 {} \\
b_1&\mapsto  a_1{} \\
a_2&\mapsto a_2b_2 a_1\\
b_2&\mapsto  a_2 \\
&\quad\vdots\\
a_\ell&\mapsto a_\ell b_\ell a_ {\ell -1} \\
b_\ell &\mapsto  a_\ell  .\\
\endaligned\right.$$ The length of $\alpha ^p(a_\ell )$ (and of its conjugacy
class) is equivalent to
$p^{\ell -1}\lambda ^p$, with
$\lambda $ the Perron-Frobenius eigenvalue of $\pmatrix
1&1\\1&0\endpmatrix$.

\subhead Polynomial growth\endsubhead

Let $n\ge2$. For the automorphism $\alpha _n$ of $F_n$ defined by $\alpha
_n(a_1)=a_1$ and
$\alpha _n(a_i)=a_ia_{i-1}$ for $2\le i\le n$, both the element $a_i$ and the
conjugacy class $\bar a_i$ grow polynomially with degree $i-1$. In
particular, $d(\hat\alpha _n)$ equals the maximal value $n-1$ (so $\hat\alpha $
is not geometric for
$n\ge3$). The rank of $\fix \alpha _n$ is 2 (it is generated by $a_1$ and
$a_2a_1a_2\mi)$.

We need a slightly different example, with every generator but one
mapped to a conjugate. We write $x^y$ for $yxy\mi$.

\nom\poly
\thm{Lemma \sta}  Let $\ell \ge1$. For  the automorphism $\beta _\ell $ of
$F_{\ell +2}=\langle a,a_0,a_1,\dots,a_\ell \rangle$ defined by :
$$\left\{\aligned a&\mapsto a\\
a_0&\mapsto a_0a\\
a_1&\mapsto a_1{}^{a_0a}\\
a_2&\mapsto a_2{}^{a_1a}\\
&\quad\vdots\\
a_\ell &\mapsto a_\ell {}^{a_{\ell -1}a},
\endaligned\right.$$
the conjugacy class of $aa_\ell $ grows with degree $\ell +1$. If one adds a
generator $t$ with $t\mapsto ta_\ell a$, the class of $t$ grows with degree
$\ell +2$.
\fthm

\example{Remark}
If an automorphism of $F_n$ maps every generator   to
a conjugate, no conjugacy class grows polynomially with degree $n-1$. This
follows from Theorem
\inal{} and Corollary \deux.
\endexample

\demo{Proof} One first shows by induction on $i\ge1$ that no cancellation
occurs when computing iterates $\beta _\ell ^p(a_i)$, because the initial
letter of
$\beta _\ell ^p(a_i)$ is $a_{i-p}$ for $p\le i-1$, and $a_0$ for $p\ge i$ (and
the final letter is the inverse of the initial one).
The length of the (non cyclically reduced
word) $\beta _\ell ^p(a_\ell )$ is the $\ell^1$-norm  of the
vector
$$\pmatrix
n_a\\n_0\\n_1\\\vdots\\
n_\ell 
\endpmatrix=
\pmatrix1&1&2&\dots&2\\
&1&2&  \\
&&\ddots&\ddots\\
&&&1&2\\
&&&&1
\endpmatrix ^p\pmatrix
0\\\vdots\\\vdots\\
0\\1
\endpmatrix.$$
It grows with degree $\ell +1$. The word $\beta _\ell ^p(aa_\ell )=a\beta
_\ell ^p(a_\ell )$ is cyclically reduced, so the class of $aa_\ell $ grows with
degree
$\ell +1$.

 The $p$-th iterate of the new generator
$t$ is a cyclically reduced word containing $a_\ell , \beta _\ell (a_\ell
),\dots,\beta _\ell ^{p-1}(a_\ell )$ as disjoint subwords, so the class of $t$
grows with degree
$\ell +2$.
\cqfd\enddemo

\subhead Mixed growth\endsubhead

For $n\ge3$, we construct an automorphism $\theta _n$ of $F_n$ with   $e$
and
$d$   as close as possible to
$(n-1)/2$. Its fixed subgroup has rank 2.

First assume $n$   odd, and write $n=2\ell +3$. Consider
$F_n=\langle a,b,a_0,a_1,b_1,\dots,a_\ell ,b_\ell \rangle$. Let $u=[a,b]$ and
$u_i=[a_i,b_i]$. Recall that $\tau :(a,b)\mapsto (aba, ba)$ is an
exponentially growing automorphism fixing $u=[a,b]$. Define $\theta
_n\in\Aut(F_n)$ by:
$$\left\{\aligned
a&\mapsto aba\\
b&\mapsto ba\\
a_0&\mapsto a_0u\\
a_1&\mapsto (a_1b_1a_1){}^{a_0u}\\
b_1&\mapsto ( b_1a_1){}^{a_0u}\\
a_2&\mapsto (a_2b_2a_2){}^{u_1u}\\
b_2&\mapsto ( b_2a_2){}^{u_1u}\\
&\quad\vdots\\
a_\ell &\mapsto (a_\ell b_\ell a_\ell ){}^{u_{\ell -1}u}\\
b_\ell &\mapsto ( b_\ell a_\ell ){}^{u_{\ell -1}u}\\
\endaligned\right.$$

Geometrically, $\theta_n$ is represented by a homotopy equivalence $\psi$ on a 2-complex $X_n$ built as follows. Take disjoint punctured tori $T, T_1, \dots,T_\ell$, with points $v,v_1,\dots, v_\ell$ on the boundary. Glue a circle to $v$, and add edges $vv_i$. The  map $\psi$ induces a pseudo-Anosov homeomorphism on each punctured torus.

The automorphism $\theta_n$ has $\ell +1$ exponential strata. Furthermore, consider the subgroup $P_0$ generated by
$u,a_0,u_1,\dots,u_\ell $ (it is the fundamental group of the 1-complex obtained from $X_n$ by removing the interior of each torus). 
It
 is $\theta _n $-invariant, and the restriction of
$\theta _n
$ to $P_0$ is
$\beta _\ell $. Thus $\hat\theta _n$ satisfies $e=d=\ell +1=\frac{n-1}2$ by
Lemma
\poly. Both inequalities  of Corollary \dcor{} are   equalities. 

If $n$ is even, we write $n=2\ell +4$ and we add a generator $t$ mapped to
$tu_\ell u$. We get an automorphism with $e=\frac n2-1$ and $d=\frac n2$ by
Lemma \poly.

The inequality $e+d\le n-1$ of Corollary \dcor{} is an equality.

These examples have only one exponential growth type $(\lambda ,0)$, with
$\lambda $ the Perron-Frobenius eigenvalue of $\pmatrix
2&1\\1&1\endpmatrix$. It is easy to modify the construction so that there
are
$e$ distinct exponential growth types $(\lambda ^ {i+1},0)$, $0\le i\le \ell$,
by using
$\tau ^{i+1}$ rather than
$\tau
$ when defining the image of
$a_i$ and $b_i$.

\subhead Optimality\endsubhead

We can now show:

\thm{Theorem \sta}  Given $(\varepsilon ,\delta )$ belonging to the closed
quadrilateral with vertices  $(0,0),(0,n-1),(\frac{n-1}2,\frac{n-1}2),
(\frac{3n-2}4,0)$, there exists $\alpha \in\Aut(F_n)$ such that:
\roster
\item Any improved relative train track map representing a power of
$\hat\alpha $ has $\varepsilon $ exponential strata. Equivalently, there are
$\varepsilon $ attracting laminations. 
\item There are $\varepsilon $ distinct exponential growth rates $(\lambda
_i,0)$ with
$\lambda _i>1$.
\item There is a conjugacy class whose growth is polynomial of degree
$\delta $.
\item The rank of $\fix\alpha $ is the maximal value $\rho _0$ permitted by
the inequalities of Corollary \pcor.
\endroster
\fthm

\demo{Proof} We write $\rho $  for $\rk\fix\alpha $.  
The inequalities of Corollary \pcor{} are
$$\aligned
\varepsilon +(\delta -1)^++\rho &\le n  \\
4\varepsilon +2\delta +2\rho  &\le 3n+1\qquad \text{($\le 3n$ if $\delta
=0$)}.
\endaligned$$
Since $(\varepsilon ,\delta )$ is in the quadrilateral, the maximal
value
$\rho _0$ is always
$\ge1$.  Which of
the two inequalities is the limiting one depends on the position of
$\varepsilon $ with respect to $n/2$ if $\delta =0$, with respect to $(n-1)/2$
if
$\delta >0$. 

We shall construct   $\alpha $ satisfying $(1)$, $(3)$, $(4)$, using the
automorphisms $\tau , \varphi _n, \alpha _n, \theta _n$ introduced above. As
in the construction of
$\theta _n$, one achieves (2) by using varying powers of $\tau $.

We write $I_\ell$ for the identity automorphism of $F_\ell$.

$\bullet$  If $\varepsilon =0$ and $\delta >0$, we define
$\alpha =\alpha _{\delta +1}*I_{n-\delta -1}$, with
$\alpha _{\delta +1}\in\Aut(F_{\delta +1})$ defined above.

$\bullet$
 Suppose $\delta =0$. The inequalities are $\varepsilon +\rho \le n$ and
$4\varepsilon +2\rho \le 3n$. 

  If
$\varepsilon \le
\frac n2$, we write $F_n$ as the free product of $F_{n-2\varepsilon }$ and
$\varepsilon $ copies of
$F_2$, and we define $\alpha =I_{n-2\varepsilon }*\tau *\dots*\tau $. Then
$\rho =n-2\varepsilon +\varepsilon =\rho _0$.

If $\varepsilon >\frac n2$, first suppose $n$ is even, say $n=2\ell $. The
required value of
$\rho $ is
$\rho _0=\frac{3n-4\varepsilon }2 =3\ell -2\varepsilon $. Write $F_n$ as the
free product of
$F_{4\varepsilon -4\ell+2}$ with $3\ell-2\varepsilon -1$ copies of $F_2$ and
define $\alpha =\varphi _{4\varepsilon -4\ell+2}*\tau *\dots*\tau $. The
number of exponential strata is
$3\varepsilon -3\ell+1+3\ell-2\varepsilon -1=\varepsilon $, and $\rho
=1+3\ell-2\varepsilon -1=\rho _0$. If $n=2l+1$, we take the free product of
$\varphi _{4\varepsilon -4\ell+1}$ with
$3\ell-2\varepsilon  $ copies of $\tau $. Then $e=3\varepsilon
-3\ell+3\ell-2\varepsilon =\varepsilon $ and $\rho =1+3\ell-2\varepsilon
=\rho _0$.

We now 
suppose
$\delta ,\varepsilon \ge1$. This implies
$n\ge3$ and $\rho _0\ge2$.

$\bullet$
Suppose $\delta $, $\varepsilon \ge1$, and $\varepsilon \le (n-1)/2$. Then 
$\rho _0=n-\varepsilon -\delta +1$. 

The construction   uses auxiliary parameters
$w,x,y,z\ge0$, to be determined in terms of $\delta $ and $\varepsilon $. 

Starting with $\theta _{2w+3}$, for which $d=e=w+1$ and $\rho =2$, we add
$x$ generators so as to obtain an automorphism of $F_{2w+3+x}$ with
$e=w+1$ and
$d=w+1+x$ (map the first generator $t_1$ to $t_1u_\ell u$ as in Lemma \poly,
then
$t_i$ to
$t_it_{i-1}$ as in the definition of $\alpha _n$). We then take the free product
with  $I_y$ and with
$z$ copies of
$\tau
$.

We get an
automorphism of a group of rank $2w+3+x+y+2z$, with 
$d=w+1+x$, 
$e= w+1+z$, and 
$\rho =2+y+z$. We must therefore find $w,x,y,z\ge0$ satisfying
$$\align w+1+x&=\delta \\
w+1+z&=\varepsilon \\
2+y+z&=\rho _0=n-\varepsilon -\delta +1
\\2w+3+x+y+2z&=n .\endalign$$

Note that the last equation is the sum of the others.
If we know
$w$, we get
$x,y,z$ by
$$\align x&=\delta -w-1\\
y&= \rho _0-\varepsilon +w-1\\
z&=\varepsilon -w-1
   .\endalign$$
We have to choose  $w$ so that $w,x,y,z$ are non-negative. This is equivalent to
$w\ge0$,
$w\ge
\varepsilon -\rho _0+1$, $w\le \delta -1$, $w\le \varepsilon -1$. Since $\delta
$ and $\varepsilon $ are $\ge1$, and
$\rho _0\ge2$, we only need to check $\varepsilon -\rho _0+1\le \delta
-1$. This holds  because $\rho
_0=n-\varepsilon -\delta +1$ and $\varepsilon \le(n-1)/2$.

$\bullet$
Finally, suppose $\delta $, $\varepsilon \ge1$, and $\varepsilon > (n-1)/2$.
Then 
$\ds\rho _0=\left[\frac{3n+1-4\varepsilon -2\delta }2\right]$.

First assume that $n$ is odd.
We   use three parameters $w,x,z$.
We first combine $\varphi _{4x+2}$ with $\theta _{2w+3}$, in the
following sense. We consider $F_{4x+2w+3}=F_{4x+2}*\langle
a_0,a_1,b_1,\dots, a_w,b_w\rangle$. Let $u$ be a generator for the fixed
subgroup of
$\varphi _{4x+2}$. We define an automorphism of $F_{4x+2w+3}$ as being
equal to $\varphi _{4x+2}$ on the first factor and mapping 
$a_0,a_1,b_1,\dots, a_w,b_w$ by the same formulas as in the definition of
$\theta _{2w+3}$. This   automorphism   has $e=3x+1+w$, $d=w+1$,
and
$\rho =2$. We then take the free product with $z$ copies of $\tau $, so as to
increase $e$ and $\rho $ by $z$.

We now have to solve:
$$\align w +1&=\delta  \\
 3x+1+w+z&=\varepsilon \\
2 +z&=\rho _0 
\\4x+2w+3+2z&=n  .\endalign$$
We have assumed $n$ to be odd.
Setting
$n=2\ell+1$, we have $\rho _0=3\ell+2 -2\varepsilon -\delta $.
It is easy to check that  $w=\delta  -1$, 
$z=\rho _0-2$,
$x=\varepsilon -\ell$ is a non-negative solution.

For $n=2\ell+2$, we   use the same construction with $\varphi _{4x+1}$ rather
than $\varphi _{4x+2}$, defining $w,z,x$ by the exact same formulas (now
$\rho _0= 3\ell+3 -2\varepsilon -\delta $,  the second
equation   is $3x+w+z=\varepsilon$, and the fourth one is $4x+2w+2+2z=n$).
\cqfd\enddemo

\bigskip

\head \sect Appendix: Growth \endhead

\subhead  More on train tracks \endsubhead

Let $f:G\to G$
be an improved
relative train track map.
We recall some more definitions
from [\BFH].

We write $\fd(\gamma )$ for the
tightened image of $\gamma $ (the reduced path homotopic  to
$f(\gamma )$ rel\. endpoints). 
A decomposition $\gamma =\gamma _1\dots\gamma _q$ is a splitting if
$\fd^p(\gamma )=\fd^p(\gamma _1)\dots\fd^p(\gamma _q)$ for all $p\ge1$
(i.e\. there is no cancellation between $\fd^p(\gamma _j)$ and $\fd^p(\gamma
_{j+1})$).  The
subpaths $\gamma _j$ are the pieces of the spitting, and we say that
$\gamma
$ splits over each $\gamma _j$. 

 If $H_i$ is
an NEG stratum, it consists of a single edge $e_i$, and $f(e_i)$ splits as
$e_i.u_i$ with $u_i$ of height $<i$. 
If $e $ is an edge in an exponential   stratum
$H_i$, then $f(e)$ has a splitting whose pieces are edges of $H_i$ or paths
of height $<i$. 

Let $H_i$ be exponential. By  aperiodicity, every edge of $H_i$ appears in
$f(e)$, for $e$ any edge of $H_i$. Any subpath of $\fd^p(e)$ is $i$-legal (for the
purposes of this appendix, this may be taken as the definition of $i$-legal). We
call
$\Delta _i$ the (finite) collection of maximal subpaths of height
$<i$ which appear in $f(e)$, for $e$ an edge of $H_i$.   If $\delta \in\Delta _i$,
then no $\fd^p(\delta )$ is a point. 
It follows from bounded cancellation that  there exists a constant $K_i$ such that, if a
path $\gamma $ contains an $i$-legal subpath with more than $K_i$ edges in
$H_i$, then $\gamma $ splits over an edge of $H_i$ (see   4.2 in [\BFH]). 

As mentioned earlier, there is a bijection between the set
of attracting laminations  of the automorphism represented by $f$ and the set of
exponential strata of
$f$.     A bi-infinite path $\gamma $ in $G$ is a leaf of the lamination associated to
$H_i$ if and only if any finite subpath is contained in some $\fd^p(e)$, with $e$ an
edge of
$H_i$.

\nom\incl
\thm{Lemma \sta} Let $f$ be an improved relative train track map. Let $H_i,H_j$
be exponential strata, with associated laminations $\Lambda _i,\Lambda _j$. 
The following are equivalent:
\roster
\item
$\Lambda _j\inc \Lambda _i$.
\item There exist edges $e_i,e_j$ of
$H_i,H_j$, and
$p\ge1$, such that
$\fd^p(e_i)$ splits over $e_j$. 
\item Given  edges $e_i,e_j$ of
$H_i,H_j$, there exists  
$p\ge1$ such that
$\fd^p(e_i)$ splits over $e_j$. 
\endroster
\fthm

\demo{Proof} $(2)\iff(3)$ follows from aperiodicity.  If $\fd^p(e_i)$ splits
over $e_j$, every $\fd^q(e_j)$ is contained in $\fd^{p+q}(e_i)$, so every leaf
of
$\Lambda _j$ is a leaf of $\Lambda _i$. Conversely, if $\Lambda _j\inc
\Lambda _i$, let $\gamma _0$ be  a   segment in a leaf of
$\Lambda _j$. It
  is contained in some $\fd ^p(e_i)$, and by   bounded
cancellation (see above)  $\fd^p(e_i)$ splits over  some
 edge of $H_j$ contained in  $\gamma _0$ if $\gamma _0$ is long enough.
\cqfd\enddemo 

\example{Remark} In particular, $j\le i$ if $\Lambda _j\inc\Lambda _i$. But
 the total order on the set of exponential strata of $f$
defined by $H_i\le H_j$ if $i\le j$ does not have an intrinsic meaning.
\endexample

 The
 Perron-Frobenius eigenvalue $\lambda _i$ of   $H_i$ is called the expansion
factor of $\Lambda _i$ (for the automorphism represented by $f$). 

Given $\Phi \in\Out(F_n)$, it is only true that some power $\Phi ^q$ is
represented by $f:G\to G$ as above. The attracting laminations of $\Phi $ are
those of $\Phi ^q$, and we define the expansion factor of $\Lambda $ for $\Phi $
as $\lambda  ^{1/q}$, where   $\lambda  $ is the expansion factor of $\Lambda $
for
$\Phi ^q$. 

\subhead Growth types\endsubhead

Given $\lambda \ge1$ and $m\in\N$, we say that a conjugacy class $\cc
g$   grows like $\lambda ^pp^m$ under $\Phi $, or has growth type $(\lambda ,m)$, 
if the length  $|\Phi ^p(\cc g)|$ grows like $\lambda ^pp^m$ in the sense that there
exist constants
$C_1,C_2>0$ with
$C_1
\lambda ^pp^m\le |\Phi ^p(\cc g)|\le C_2 \lambda ^pp^m$ for all $p\ge1$.   The
set of growth types is ordered lexicographically, so that $(\lambda
,m)\le(\lambda ',m')$ if $\lambda ^pp^m\le \lambda '{}^pp^{m'}$ as $p\to \infty$. 

We define similarly the growth type of an element
$g
\in F_n$ under
$\alpha \in\Aut(F_n)$, and of  an edge-path (or a loop) $\gamma $ in $G$ under a
relative train track map $f:G\to G$ by considering the simplicial length
$|\fd^p(\gamma )|$. 

If $g$ or $\cc g$ grows like
$\lambda ^pp^m$ under some positive power $\alpha ^q$ or $\Phi ^q$, it grows
like $(\lambda ^{1/q})^pp^m$ under $\alpha $ or $\Phi $. This   allows us to
replace an automorphism by a power whenever convenient.

We consider the set of attracting laminations of $\Phi $, ordered by inclusion.
Each attracting lamination $\Lambda $ has an expansion factor $\lambda _0>1$
(see above). 
 From this data, we shall now
associate to each $\Lambda $ a growth type $c=(\lambda ,m)$, with $\lambda
>1$ and $m\in \N$.
The definition will ensure that, if
$\Lambda $ is associated to an exponential stratum $H_i$ of $f:G\to G$
representing $\Phi $, and
$e$ is an edge of $H_i$, then the length of $\fd^p(e)$     grows like
$\lambda ^pp^m$ (see Proposition 6.4). 

The definition is by induction on the number of laminations contained in
$\Lambda $, using the following rules.  Let $\lambda
_0$ be the expansion factor of $\Lambda $. If
$\Lambda
$ is minimal (for inclusion), then
$c=(\lambda _0,0)$.  If not,
let
$(\lambda ',m')$ be the maximum growth type for $\Lambda '\subsetneq \Lambda
$. If
$\lambda '<\lambda _0$, then $c=(\lambda _0,0)$. If $\lambda '>\lambda _0$, then
$c=(\lambda ',m')$. If $\lambda '=\lambda _0$, then $c=(\lambda _0,m'+1)$.  
In all cases,
$c(\Lambda )=\max_{\Lambda '\subseteq\Lambda }c(\Lambda ')$. Also note that,
if $m>0$, there is a decreasing chain of laminations $\Lambda =\Lambda
_0\supsetneq \Lambda _1\supsetneq \dots\supsetneq \Lambda _m$, so that
$m\le s$.

\example{Example} Consider the automorphism of $F_4$ defined by $a\mapsto 
abaa'$,
$b\mapsto ba$,
 $a'\mapsto  a'b' $, $b'\mapsto a'$. There are two attracting laminations $\Lambda
$ and $\Lambda '$, with $\Lambda '\inc\Lambda $. The expansion factor $\lambda_0
$ of
$\Lambda $ is the Perron-Frobenius eigenvalue $\mu $ of $\pmatrix
2&1\\1&1\endpmatrix$, the expansion factor $\lambda '_0$ of
$\Lambda '$ is the eigenvalue $  \nu  $ of $\pmatrix
1&1\\1&0\endpmatrix$. One has $\mu >\nu $; the growth type
of
$\Lambda
$ is $(\mu ,0)$, corresponding to the fact that $|\Phi ^p(\cc a)|$ grows like $\mu
^p$.  For   the automorphism   $a\mapsto  aba'$, $b\mapsto
a$,
 $a'\mapsto  a'b'a' $, $b'\mapsto a'b'$, one still has two laminations
$\Lambda '\inc\Lambda $, but now $\lambda _0=\nu
$ and 
$\lambda '_0=\mu
$, so $\lambda _0<\lambda '_0$. The growth type of $\Lambda $ is $(\mu ,0)$,
and $|\Phi ^p(\cc a)|$ grows like $\mu ^p$. Finally, consider $a\mapsto  aba'$,
$b\mapsto a$,
 $a'\mapsto  a'b' $, $b'\mapsto a'$. One has $\lambda _0=\lambda '_0=\nu $. The
growth type of $\Lambda $ is $(\nu ,1)$, and $|\Phi ^p(\cc a)|$ grows like $p\nu ^p$.
\endexample

\nom\gro
\thm{Theorem \sta} Let $\Phi \in\Out(F_n)$. \roster
\item Given $g\in F_n$, there exist $\lambda \ge1$ and $m\in\N$ such that
$|\Phi ^p(\cc g)|$ grows like $\lambda ^pp^m$.  
\item If $|\Phi ^p(\cc g)|$ grows like $\lambda ^pp^m$, and $m'<m$, there exists
$g'$ such that $|\Phi ^p(\cc g')|$ grows like $\lambda ^pp^{m'}$.
\item Given $(\lambda ,m)$ with $\lambda >1$, there exists $g$ such
that
$|\Phi ^p(\cc g)|$ grows like $\lambda ^pp^m$ if and only if $(\lambda ,m)$ is
the growth type of some attracting lamination of $\Phi $.
\endroster
\fthm

When $\lambda >1$, the power $m$ which appears in the growth of $|\Phi
^p(\cc g)|$ satisfies
$m\le s \le  n/2-1$ by Theorem \cha. When $\lambda =1$, the maximum value of
$m$   is the number that we have called $d$. It satisfies $d\le n-1$ (see Corollary
\dcor).

Before proving Theorem \gro, we also note:

\thm{Corollary \sta} Given $\alpha \in\Aut(F_n)$ and $g\in F_n$, there exist
$\lambda \ge1$ and $m\in \N$ such that $|\alpha ^p(g)|$ grows like $\lambda
^pp^m$.
\fthm 

\demo{Proof} Extend $\alpha $ to an automorphism $\beta  $ of $F_{n+1}$ by
mapping the new generator $t$ to itself. The growth of $g$ under $\alpha $ is
that of the conjugacy class $\cc {tg }$ under $\hat\beta $. \cqfd\enddemo

 The set of growth types of elements of $F_n$ under $\alpha $ is the same as
the set of growth types of conjugacy classes under $\hat \alpha $, except
that there may be elements with growth $p^{d+1}$ (see Lemma \trois).

\subhead Proof of Theorem \gro\endsubhead
 
The rest of this section is devoted to  the proof of Theorem \gro. See [\BFHd,
\BH, \LLD] for partial results. Our proof elaborates on an argument due  to
Bridson-Groves.

After replacing $\Phi $ by a power, we may assume that it is represented by
  an improved relative train track map $f:G\to G$. The heart of the
proof  is to show that {\it any edge-path or loop $\gamma
$ in
$G$  has a growth    type $(\lambda ,m)$},  
  in the sense that
 the length of
$\fd^p(\gamma )$ is bounded between $C_1\lambda ^pp^m$ and $C_2\lambda
^pp^m$ for some $C_1,C_2>0$. 

Recall that we have defined a growth type $ c=(\lambda ,m)$   for an attracting
lamination
$\Lambda $, hence also for an exponential stratum $H_i$. We write $c_i$ for the
growth type attached to $H_i$. 

\nom\grop
\thm{Proposition \sta}  Given   $\gamma $, let 
$C_\gamma $ be the set  of all
$j$  such that some $\fd^p(\gamma )$ splits over  an edge belonging to an
exponential stratum $H_j$.   If $C_\gamma =\ev$, then  $\gamma $ grows like
$p^m$ for some $m\in\N$. Otherwise, the growth type of $\gamma $ is
the maximal $c_j$, for $j\in C_\gamma $. \fthm

In particular, the growth type of an
edge in an exponential stratum
$H_i$  is
$c_i$ by Lemma \incl{} and the equality $c(\Lambda )=\max_{\Lambda
'\subseteq\Lambda }c(\Lambda ')$.

The proof of the proposition is by induction on the height of 
$\gamma
$.
First suppose that $\gamma $ is a single edge $e$ in an
exponential stratum $H_i$, and the proposition is true for paths of height $<i$.
Let
$c_i=(\lambda , m)$, and let
$\lambda _i\le\lambda $ be the Perron-Frobenius eigenvalue attached to $H_i$. 
We show that $|\fd^p(e)|$ grows   like $\lambda ^pp^m$.

As mentioned above,  $f(e)$ splits over edges of $H_i$ and paths $\delta
\in\Delta _i$  (recall that we have defined $\Delta _i$ as the finite set of maximal
subpaths of height
$<i$ which appear in $f(e')$ for $e'$ an edge of $H_i$).  Thus $\fd^p(e)$ splits
over edges of $H_i$ and paths of the form $\fd^{q}(\delta )$ with $q\le
p-1$ and $\delta \in\Delta _i$. Up to constants (which we will not write), the
number of edges in
$\fd^p(e)\cap H_i$ is
$ \lambda _i^p$, and for given $q\le p-1$ and $\delta \in\Delta _i$ the number of
subpaths
$\fd^{q}(\delta )$ is
$\lambda _i^{p-q-1}$ (as they are created by edges in $\fd^{p-q-1}(e)\cap H_i$).

We first show that the length of $\fd^p(e)$ grows at most like $\lambda ^pp^m$.
It suffices to show that, for a given $\delta $ in  the   set $\Delta _i$, the
sum
$\sum_{q=1}^p \lambda _i^{p-q}|\fd^{q}(\delta )|$ grows at most
like
$\lambda ^pp^m$. If some $\fd^{q}(\delta )$ splits over an edge in an exponential
stratum
$H_j$, then so does some $\fd^{q'}(e)$. By the induction hypothesis,
$|\fd^{q}(\delta )|$ grows either polynomially or with growth type $c_j=(\lambda
_j,m_j)$, with
$j\in C_\gamma $ and $j<i$, so we have to show that $S_p= \sum_{q=1}^p \lambda _i^{p-q}\lambda _j^qq^{m_j}$ grows at most like $\lambda ^pp^m$. 
But the growth type of a lamination was defined in such a
way that this holds, since  $S_p$ grows like $\lambda_i^p $ if $\lambda_j<\lambda_i$, like $\lambda _j^pp^{m_j}$ if $\lambda_j>\lambda_i$,   like $\lambda _j^pp^{m_j+1}$ if $\lambda_j=\lambda_i$.

We now show that $e$ grows
  at least  like $\lambda ^pp^ m$. 
This is clear if $\Lambda _i$ is minimal (for
inclusion). If $\lambda  >\lambda _i$, there is $\Lambda
_j\varsubsetneq\Lambda _i$ with $c_j=c_i$. By Lemma \incl, some $\fd^q(e )$
splits over an edge $e_j$ of $H_j$, so $e $ grows as fast as $e_j$, whose growth
type is
$c_j$ by induction. The  case   when $\lambda  =\lambda _i$ is harder.

The result is clear if $m=0$, so assume $m>0$. 
Then there is $\Lambda _j\varsubsetneq \Lambda _i$ with 
$c_j=(\lambda ,m-1)$. Some $\fd^q(e )$ splits over an edge $e_j$ of $H_j$. This
edge is contained in some
$\fd^{q'}(\delta )$, with
$q'\le q-1$ and
$\delta \in
\Delta _i$. 
The splitting of $\fd^q(e )$
over $e_j$ may not be compatible with the splitting of $\fd^q(e )$ over
$\fd^{q'}(\delta )$. 
But further iterates of $e$ split over long $j$-legal paths,  and by
bounded cancellation we may therefore assume (after increasing
$q$)  that a subpath  $\delta '=\fd^{q'}(\delta )$  of $\fd^q(e
)$ splits over an edge of $H_j$, hence has growth type at least $(\lambda
,m-1)$ by induction. 

Fix $q$ and $\delta '$, and  consider $\fd^{qp}(e
)$. For each $\ell<p$ it contains (up to constants) at least $\lambda ^{q(p-\ell)}$
subpaths
$\fd^{q\ell}(\delta ')$, each with length   $\lambda ^{q\ell}({q\ell})^{m-1}$. It
follows that the length of $\fd^{qp}(e
)$ is bounded below by   $\sum_{\ell=1}^{p-1}\lambda
^{q(p-\ell)}\lambda ^{q\ell}({q\ell})^{m-1}$, hence by $\lambda ^{qp}(qp)^m$.

This completes the proof of the  induction step for $\gamma $   an edge in an
exponential stratum. In general, we now know that $ \gamma $ grows at least like
the maximal $c_j$. For the upper bound, we use the following fact.

\thm{Lemma \sta{} [\BFH, Lemmas 4.1.4, 4.2.6, 5.5.1]} If $\gamma $ has height $i$,
there exists
$p_0$ such that 
$\fd^{p_0}(\gamma )$ has a splitting whose pieces are
edges of $H_i$, paths of height $<i$, Nielsen paths, and exceptional paths. \cqfd
\fthm

  Nielsen paths do not grow, and exceptional paths grow linearly. Paths of height
$<i$ grow at most like the maximal $c_j$ by the induction hypothesis. We
know how edges in an exponential stratum grow, so there only remains the case
when
$\gamma $ is an edge   in an NEG stratum. 
We have
$\fd^p(\gamma )=\gamma u\fd(u)\fd^2(u)\dots\fd^{p-1}(u)$ for some path $u$ of
height
$<i$.  By
the induction hypothesis, 
$u$ grows
either polynomially or with growth type $c_j$ with $j\in
C_\gamma $.
Thus $\gamma $ grows at the same speed as $u$ if $u$ grows exponentially,
like
$ p^{m+1}$ if $u$ grows like $p^m$.

The proof of Proposition \grop{} is now complete.  To prove the first
assertion of Theorem \gro, we simply observe that  the type of growth of a
conjugacy class $\cc g$ (under $\Phi
$) is the same as that of the loop representing $\cc g$ in $G$  (under $f$). If $H_i$
is an exponential stratum,  any conjugacy class represented by an $i$-legal
loop meeting $H_i$ has growth type  
$c_i$. This shows assertion (3), and assertion (2) when $\lambda >1$. Assertion
(2) in the polynomial case is well-known, but we sketch a proof for
completeness. If
$\cc g$ grows like $p^m$, the arguments given above imply the existence of    a
path
$u$ which grows like
$p^{m-1}$. This path is a loop by [\BFH], but images $\fd^p(u)$ may fail to be
reduced as loops. It is easy to check, however,  that for each $p$ one of the loops
$\fd^p(u)$, $\fd^{p+1}(u)$, $\fd^p(u) \fd^{p+1}(u)$ is reduced, and assertion (2)
follows.

\bigskip\null
\Refs 
\widestnumber\no{99}
\refno=0

\bref \by M. Bestvina, M. Feighn, M. Handel\paper The Tits alternative for
$\Out(F_n)$,  I, Dynamics of exponentially-growing
       automorphisms \jour Ann. of Math. \vol 151 \yr2000\pages 517--623
\endref

\bref \by M. Bestvina, M. Feighn, M. Handel\paper The Tits alternative for
$\Out(F_n)$,  II, A Kolchin type theorem \jour Ann. of Math. \vol 161
\yr2005\pages 1--59
\endref

\bref \by M. Bridson, D. Groves\paper
The quadratic isoperimetric inequality for mapping tori of free group automorphisms
\jour Memoirs of the AMS (to appear)
\endref

\bref \by M. Bestvina, M. Handel\paper Train tracks for automorphisms
of the free group \jour Ann. Math.\vol135 \yr1992\pages1--51 \endref

\bref  \by D.J. Collins, E.C. Turner\paper All automorphisms of free groups
with maximal rank fixed subgroups\jour Math.
Proc. Camb. Phil. Soc. \vol119\yr1996\pages615--630\endref

 \bref \by  J.L. Dyer, P.G. Scott \paper
Periodic automorphisms of free groups \jour
Comm. Algebra \vol 3 \yr 1975 \pages 195--201  
\endref

\bref\by D. Gaboriau, A. Jaeger, G. Levitt, M. Lustig\paper An index for counting
fixed points of automorphisms of free groups\jour Duke Math.
Jour.\vol93\yr1998\pages425--452 
\endref

\bref  \by D. Gaboriau, G. Levitt\paper The rank of actions on \Rt s
\jour Ann. Sc. ENS \vol28\yr1995\pages549--570 \endref

\bref   \by D. Gaboriau, G. Levitt, M. Lustig\paper A dendrological
proof of the Scott conjecture for automorphisms of free groups
\jour    Proc. Edinburgh Math. Soc. \vol  41  \yr1998 \pages 325--332 
 \endref

\bref \by F. Gautero, M. Lustig
\paper
  The mapping-torus of a free group
automorphism is hyperbolic relative to the canonical subgroups of polynomial
growth 	\jour arXiv:0707.0822
\endref

\bref\by G. Levitt, M. Lustig\jour Unpublished notes, 1998
\endref

\bref\by G. Levitt, M. Lustig\paper Periodic ends, growth rates, H\"older
dynamics for automorphisms of free groups \jour Comm. Math.
Helv.\yr2000\vol75\pages415--430\endref

\bref \by A. Piggott  \paper  Detecting the growth of free group
automorphisms by their action on the homology of subgroups of finite index
\jour math.GR/0409319
\endref

\bref  \by H. Short\paper Quasiconvexity and a theorem of Howson's
\inbook Group theory from a geometrical viewpoint\eds Ghys, Haefliger,
Verjovsky\publaddr World Scientific\yr 1991
\pages 168--176\endref

\endRefs

\medskip

\address  Gilbert Levitt: LMNO, umr cnrs 6139, BP 5186, Universit\'e de Caen,
14032 Caen Cedex, France.\endaddress\email
levitt\@math.unicaen.fr{}{}{}{}{}\endemail

\enddocument